\documentclass[12pt]{amsart}
\usepackage{amssymb}
\usepackage{latexsym}
\usepackage{epsf}
\newtheorem{thm}{Theorem}[section] 
\newtheorem{pro}[thm]{Proposition}
\newtheorem{lem}[thm]{Lemma}
\newtheorem{cor}[thm]{Corollary}

\theoremstyle{definition}
\newtheorem{dfn}{Definition}[section]

\def\1{{\rm1\mathchoice{\kern-0.25em}{\kern-0.25em}
        {\kern-0.2em}{\kern-0.2em}I}}

\newcommand{\one}{{ \rm \setlength{\unitlength}{1em}
\begin{picture}(0.75,1)
\put(0,0){$1$}\put(0.34,0){\line(0,1){0.65}}
\end{picture} }}
\newcommand{\oneun}{\one\!{}_1}
\newcommand{\lmn}[1]{\vadjust{\setbox1=\vtop{\hsize 25mm
\parindent=0pt\baselineskip=9pt
\rightskip=4mm plus 4mm#1}
\hbox{\kern-26mm\smash{\raise .5ex\box1}}}}

\newcommand{\psdiag}[3]{\hspace{1mm}\raisebox{-#1mm}{\epsfysize#2mm
\epsffile{#3.eps}}\hspace{1mm}}

\newcommand{\nc}{\newcommand}
\def\be#1\ee{\begin{equation}#1\end{equation}}
\nc{\bc}{\begin{center}}
\nc{\ec}{\end{center}}
\nc{\bb}{\mathbb}
\nc{\cal}{\mathcal}
\nc{\N}{{\mathsf N}}
\nc{\K}{{\mathsf K}}
\nc{\fk}{\mathbf{k}}
\theoremstyle{remark}

\def\v8{\vskip 8pt}
\def\a{\alpha}

\def\la{\langle}
\def\ra{\rangle}
\def\l{\lambda}
\def\n{\nu}
\def\m{\mu}

\def\S{\mathcal S}
\begin{document}

\title[BCD modular categories]{Modular categories
 of types B,C and D}
\author{Anna Beliakova}
\address{Mathematisches Institut, Universit\"at Basel, Rheinsprung 21,
CH-4051 Basel, Switzerland}
\email{Anna.Beliakova@unibas.ch}
\author{Christian Blanchet}
\address{L.M.A.M., Universit\'e de Bretagne-Sud,
Centre de Recherches Tohannic, BP 573, F-56017 Vannes, France  }
\email{Christian.Blanchet@univ-ubs.fr}

\begin{abstract}
We construct four series of modular categories
from the two-variable Kauffman polynomial, without use 
of  the representation theory of
quantum groups at roots of unity.
The specializations of this polynomial corresponding to
quantum groups of types B, C and D  produce series of
pre-modular categories. One of them 
turns out to be  modular and three others
satisfy  Brugui\`eres' modularization criterion.
For these four series  we compute the Verlinde formulas,
and discuss spin and cohomological refinements.
\end{abstract}

\maketitle
\section*{Introduction}

Modular categories are  tensor categories  with additional 
 structure
(braiding, twist, duality, a finite set of dominating simple objects 
satisfying a
non-degeneracy axiom). 
If we remove the last axiom, we get a pre-modular category.
A pre-modular category
   provides   invariants of links, tangles, and 
sometimes of $3$-manifolds. Any modular category yields
a Topological Quantum Field Theory (TQFT) in dimension three
\cite{Tu}.

In this paper we give an elementary
 construction of  modular and pre-modular categories 
arising from the Kauffman skein relations, without use of
 the representation theory of quantum groups at 
roots of unity. Our method is  based on the
skein-theoretical construction of idempotents in  the 
Birman-Murakami-Wenzl (BMW) algebras  given in \cite{BB}.
 This work follows  the  program
of Turaev and Wenzl \cite{TW,TW2}.
We give four specifications of parameters $\alpha$ and $s$ 
(entering  the Kauffman skein relations) which lead to different 
series of modular
categories.
 In each case,
the quantum parameter $s$ is a root of unity and $\pm \alpha$
is a power of $s$. The order $l$ of $s^2$ plays a key role in the discussion.
 When $l$ is odd, then either $s^l=-1$ or $s^l=1$. We note that the two
cases are quite different: only one of them lead to a modular category, the
other one produces a non-modularizable pre-modular category.

 It is well-known that the link invariant associated
with  the fundamental representation of the quantum group of type $A_n$ 
is a specialization of the Homfly
polynomial. Taking the fundamental representations
of the quantum groups of types $B_n$, $C_n$ or $D_n$
one obtains   specializations of the Kauffman polynomial \cite{Tu1}.
 More generally, with each of these quantum groups
 at  a root 
of unity $q$ a 
pre-modular category can be associated  \cite{Kir}.
The order of $q$ determines the {\it level} $k$ of the category. 
It turns out that  categories obtained from the quantum groups
of types $A_n$ and $A_k$,
where $q$ is $(n+k)$th root of unity, are isomorphic;
 here one has to consider either a non standard
choice of the framing parameter, or the {\em projective }
subcategory.
 The isomorphism interchanges the rank $n$ and the level $k$
of the category and it is known as the level-rank duality.
 This duality has no natural explanation in the context of quantum
groups, because the roles of the parameters $n$ and $k$  
are completely different there. 

In our setting, both parameters $n$ and $k$ serve to restrict
the size of the Weyl alcove, and we have  natural 
symmetries interchanging them. Therefore, each of our
(pre-) modular categories has its level-rank duality partner.
In fact, all our specializations of parameters
can be interpreted in two different ways as a quantum group specialization.
Accordingly, we denote our categories by pairs of the 
letters $B$, $C$ and $D$ (we use just one of them if both coincide).
Our main results can be formulated as follows.

\begin{itemize}
\item We recover the symplectic 
($C$  in our notation)
 and $BC$ series of modular categories
already obtained by Turaev and Wenzl \cite{TW2}.
 These series are constructed by killing negligible
morphisms in the idempotent completed
Kauffman category. In the $BC$ case we further use
Brugui\`eres' modularization procedure \cite{Bru}.
 This could be avoided here by considering
a subcategory (see \cite[9.9]{TW2}).
\item
We obtain two new series of modular categories in the orthogonal case:
 one in the
even orthogonal case ($D$ series)
 and one in the mixed  odd-even orthogonal
case ($BD$ series).
 All of them are constructed
by using Brugui\`eres' modularization procedure.
\item
Except for the even orthogonal categories, 
we describe  explicitly
 the representative sets of simple objects and
  state the Verlinde formulas,
which give the dimensions of the  TQFT modules. 
 In the even orthogonal case, the complete description of the set 
of simple objects
depends on a tricky computation which has still to be done.

\item 
We find a correspondence
between our categories and categories obtained
by the quantum group method. We show that 
the categories constructed here give a complete set of 
invariants that can be obtained  from
quantum groups of types $B$, $C$ and $D$
by using non-spin
modules.

\end{itemize}


The paper is organized as follows.
In the first section we give the general definitions
and theorems concerning  pre-modular and modular categories.
 This includes Brugui\`eres' modularization criterion,
and an explicit description of
a modularization functor for a modularizable pre-modular
category whose {\em transparent} simple objects
are invertible.
 In the second section we recall the main definitions and  properties of 
the minimal 
 idempotents in the BMW algebras constructed in \cite{BB}. In the third 
 section we construct the completed BMW category and use it 
in order to
define series of pre-modular categories.
 In Section 4, studying    transparent  objects 
in these categories, we
 show that the symplectic category is modular and 
three  other series
 satisfy Brugui\`eres' modularization criterion.
 Then for modular categories we describe the representative 
sets of simple objects,
  give the Verlinde formulas and discuss spin and 
cohomological refinements.
In the last section we explain how our pre-modular categories
can be interpreted in terms of
quantum groups.

{\bf Conventions.}  The  manifolds throughout this paper are 
 compact, smooth and oriented. 
By a {\it  link}
we mean an { isotopy class of  an unoriented framed link}.
 Here, a framing is a non-singular normal vector field, up
to homotopy. By a {\it  tangle} in a $3$-manifold $M$
we mean {an  isotopy class of a framed tangle relative to the  boundary}.
Here the boundary of the tangle is a finite set
of points in $\partial M$, together with a nonzero vector tangent to 
 $\partial M$ at each point.
 Note that a framing together with an orientation is equivalent
to a trivialization of the normal bundle, up to homotopy.
By an {\it oriented link} we mean an isotopy class of a link 
together with a trivialization of the normal bundle, up to
homotopy.
By an {\it oriented tangle } we mean an isotopy class of a tangle 
  together with a trivialization of the normal bundle, up to
homotopy
 relative to the boundary. 
Here the boundary of the tangle is a finite set
of points in $\partial M$, together with a trivialization of the tangent
space  to 
 $\partial M$ at each point.
 In the figures, a convention using the plane
gives the  preferred framing (blackboard framing).

{\bf Acknowledgments.} The first and
the second authors thank, respectively,  the
Laboratoire de Math\'ematiques  de l'Universit\'e de Bretagne-Sud in Vannes
and 
the Mathematisches Institut der Universit\"at Basel
for their hospitality.  The authors also wish to thank Alain Brugui\`eres 
and Thang Le  explaining  algebraic
structures relevant to their constructions and for 
useful remarks on the preliminary
version of this paper.
 This work was supported in part by the Sonderprogramm zur F\"orderung
des akademischen Nachwuchses der Universit\"at Basel.


\section{Pre-modular categories and modularization}\label{one}
\subsection{Pre-modular and modular categories}
 A ribbon category is a category equipped with a tensor product,
braiding, twist and duality satisfying compatibility conditions
\cite{Tu}.
 If we are given a ribbon category ${A}$, then we can define
an invariant of links whose components are colored by  objects of
${A}$. This invariant extends to a representation of the
${A}$-colored tangle category and more generally
to a representation of the category of ${A}$-colored ribbon graphs
\cite[I.2.5]{Tu}.
 Using the ribbon structure of $A$, we get traces of morphisms and dimensions 
of objects, for which we will use the terminology {\it
quantum trace} and {\it quantum dimension}. More precisely,
for any $X\in Ob(A)$ and $f\in End(X)$
we denote by $\la f\ra \in End(\mathrm{trivial\ object})$
 the quantum trace of $f$ and by 
$\la X\ra=\la \one\!_X\ra $ the quantum dimension of $X$.
Throughout this paper
  $\one\!_X$ denotes the identity morphism of $X$.

Let $\fk$ be a field. A ribbon category will be said to be
$\fk$-linear if the Hom sets are 
$\fk$-vector spaces, composition and tensor product are
bilinear, and $End(\text{trivial object})=\fk$.
We call an object $X$ of $A$ simple if
 the map  $u\mapsto u \one\!_{X}$
from $\fk=End(\mathrm{trivial\ object})$ to $End(X)$ is an isomorphism.

\begin{dfn}
A {\em modular category} \cite{Tu}, over the field $\fk$,
is a $\fk$-linear ribbon category in which there
exists a finite family $\Gamma$ of simple objects $\lambda$
satisfying the four axioms below.
\begin{enumerate}
\item (Normalization axiom) The trivial object 
is in $\Gamma$.
\item (Duality axiom) For any object $\lambda\in \Gamma$,
its  dual  $\lambda^\ast$ is isomorphic to an object in
$\Gamma$.
\item (Domination axiom) For any object $X$ of the category
there exists a finite decomposition
$\one\!_X=\sum_i f_i\one\!_{\lambda_i}g_i$,
with $\lambda_i\in\Gamma$ for every $i$.
\item (Non-degeneracy axiom) The following matrix is invertible.
$$S=(S_{\lambda\mu})_{\lambda,\mu\in \Gamma}\;\, ,$$
where $S_{\lambda\mu}\in \fk$ is the endomorphism
of the trivial object
associated with the $(\lambda,\mu)$-colored,
$0$-framed Hopf link  with linking $+1$. 
\end{enumerate}
\end{dfn}

It follows that $\Gamma$ is a representative
set of isomorphism classes of simple objects.
 If we remove the last axiom, we get a definition of 
a {\it pre-modular} category.

\begin{dfn}
An  object $\l$ of a pre-modular category $A$
 is called {\em transparent},
if 
for any object $\mu$ in $A$
$$\psdiag{4}{12}{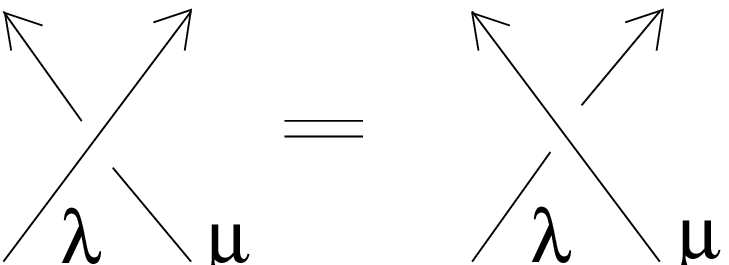}.$$
\end{dfn}
\noindent 
Such an object is  also called a central object.
 It is enough to have the above equality
for any $\mu$ in a representative set of simple objects.
 Note that a category containing a nontrivial transparent 
simple object can not be modular,
simply because the row in the $S$-matrix corresponding to this transpa\-rent
object is colinear to the row of the trivial one.
In the next subsection we show that  the absence of nontrivial 
transparent simple objects implies (under a mild assumption)
that the category is modular.

\subsection{Properties of  pre-modular categories}\label{onetwo}
We will first give some general facts about pre-modular categories.
 Let  $A$ be 
a  pre-modular category
and let $\Gamma(A)$ be a representative set of isomorphism
classes of its simple objects. We denote by
$\omega$
the {Kirby color}, i.e.
$\omega=\sum_{\l\in\Gamma(A)}\la \l \ra \l$.
 We use here the same notation as before for traces and dimensions.
In addition,
 we suppose that $A$ has no nontrivial negligible morphisms
(we quotient out by negligible morphisms if necessary). 
Note that a morphism $f\in Hom_A(X,Y)$ is called negligible if
for any $g\in Hom_A(Y,X)$ $\la fg\ra=0$.

\begin{pro}\label{slidprop} (Sliding property)
The following two morphisms in $A$ are equal.
$$\psdiag{10}{33}{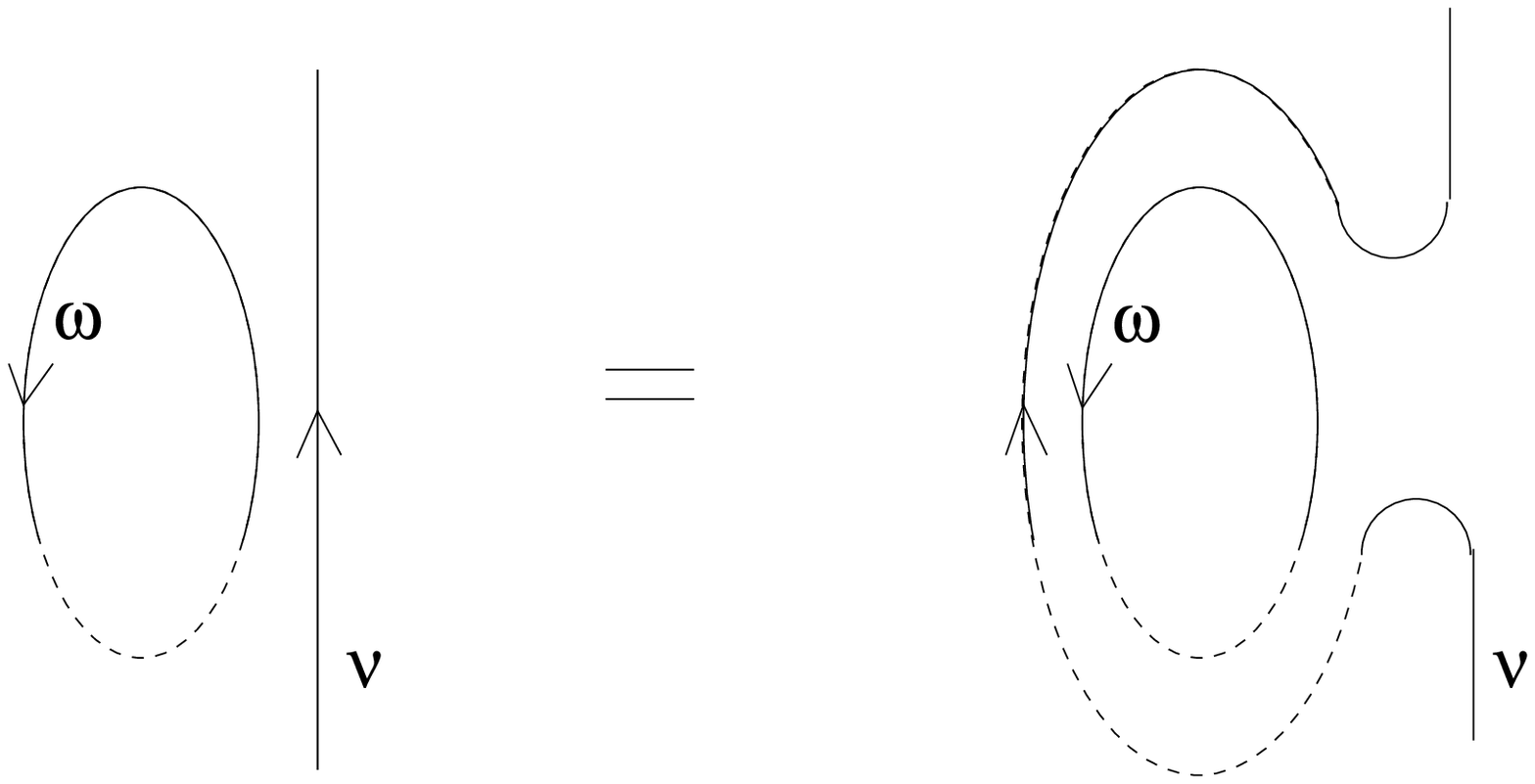}$$
\end{pro}
\noindent
Here the  dashed line represents a part of  the closed component 
colored by $\omega$. This part
can be knotted or linked with other components of 
a ribbon graph representing the morphism.
 Note that the morphism is unchanged if we reverse
the orientation of this closed component.
\begin{proof}
For $c_i, d_j\in {\Gamma (A)}$, $i=1,...,n$, $j=1,...,m$,
 we put
$$Hom_A(c_1\otimes...\otimes c_n,d_1\otimes...\otimes d_m):=
H^{d_1...d_m}_{c_1...c_n} .$$ 
With this notation
the modules $H^{\l\n}_\m$, $H^\l_{\m\n^\ast}$, $H^{\m^\ast\l}_{\n^\ast}$,
$H^{\m^\ast}_{\n^\ast\l^\ast}$, $H^{\n\m^\ast}_{\l^\ast}$ and
 $H^\n_{\l^\ast \m}$
 are mutually isomorphic, as well as the modules
$H_{\m\n^\ast\l^\ast}$, $H^{\l\n\m^\ast}$ and all obtained from them
by cyclic permutation of colors.
For example, the map $\Psi: H^{\l\n}_\m\to H^\l_{\m\n^\ast}$ and its inverse
are depicted below.
\vspace*{0.2cm}
$$\psdiag{6}{18}{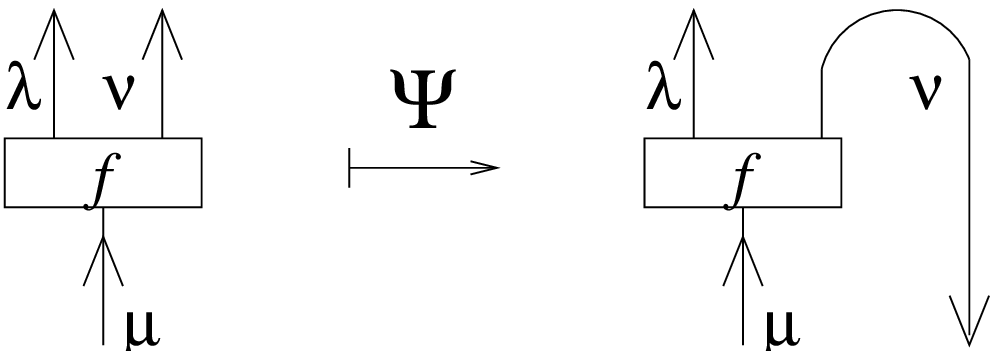}\;\;;\;\;\;\;\;\;\;\psdiag{6.5}{19.5}{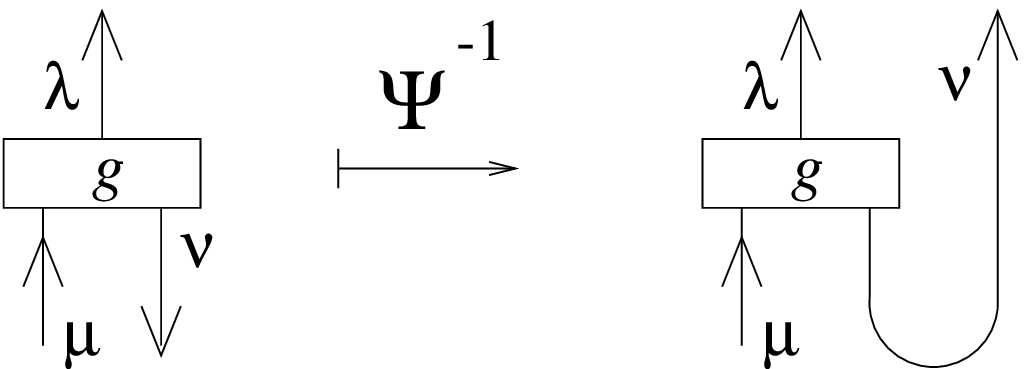}$$
\vspace*{0.2cm}

\noindent
Identifying these modules along the isomorphisms we get a symmetrized
multiplicity module ${\tilde H}^{\l\n\m^\ast}$; here only the cyclic order
of colors is important. We will represent the elements of 
${\tilde H}^{\l\n\m^\ast}$
by a circle with one incoming line (colored with $\m$) and two outgoing ones
(colored with $\l$ and $\n$), the cyclic order of lines is $(\l\n\m)$.
The module ${\tilde H}^{\m\n^\ast\l^\ast}$ is dual to
${\tilde H}^{\l\n\m^\ast}$. The natural pairing is non-degenerate,
since we have no negligible morphisms.
 We denote by $a_i$, $i\in I^{\l\n\m^\ast}$, a basis
of ${\tilde H}^{\l\n\m^\ast}$, 
and by $b_i$ the dual basis with respect to this pairing.
 Applying the domination axiom we get
that the natural map 
$\oplus_\mu \tilde H^{\l\n\m^\ast}
\otimes \tilde H^{\m\n^\ast\l^\ast }\rightarrow H^{\l\n}_{\l\n}$
is an isomorphism.  By writing the identity of 
$\l\otimes\n$ in the basis corresponding to
 $(a_i\otimes b_j)$, we get  the following decomposition formula
({\em fusion } formula)

\be \label{tprod}\sum_{\m}\sum_{i\in I^{\l\n\m^\ast}}\la\m\ra
\ \ \psdiag{9}{27}{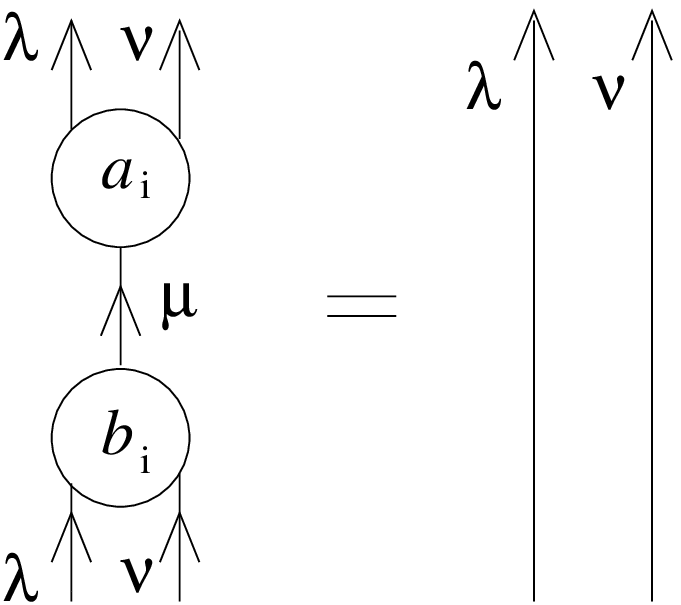}\ .\ee

The calculations below  establish the sliding property.

\begin{eqnarray*}\sum_{\l\in\Gamma (A)}\la \l\ra \psdiag{8}{24}{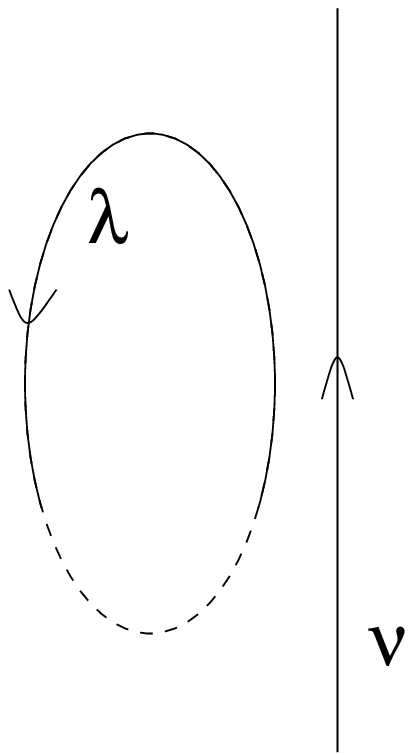}&=&
\sum_{\l,
\m}\sum_{i\in I^{\l\n\m^\ast}}\la\l\ra \la\m\ra
\psdiag{8}{24}{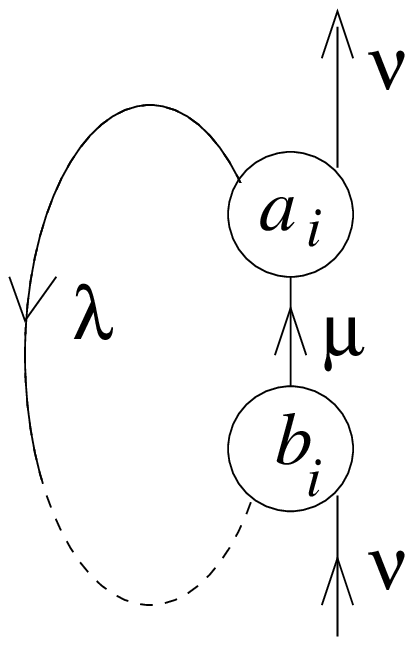}\\
&=&
\sum_{\l ,\m}\sum_{i\in I^{\l\n\m^\ast}}{\la\l\ra}{\la\m\ra}
 \psdiag{10}{30}{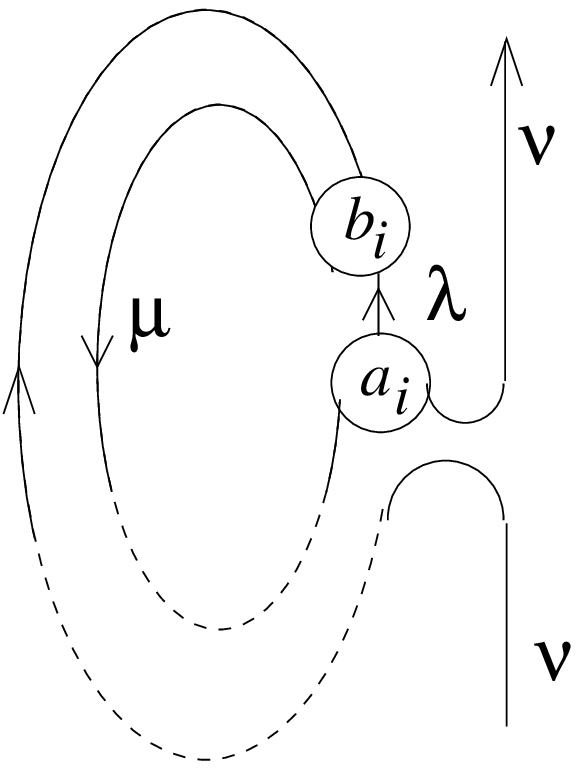}\\
&=&
\sum_{\m\in {\Gamma (A)}}\la\m\ra\psdiag{8}{24}{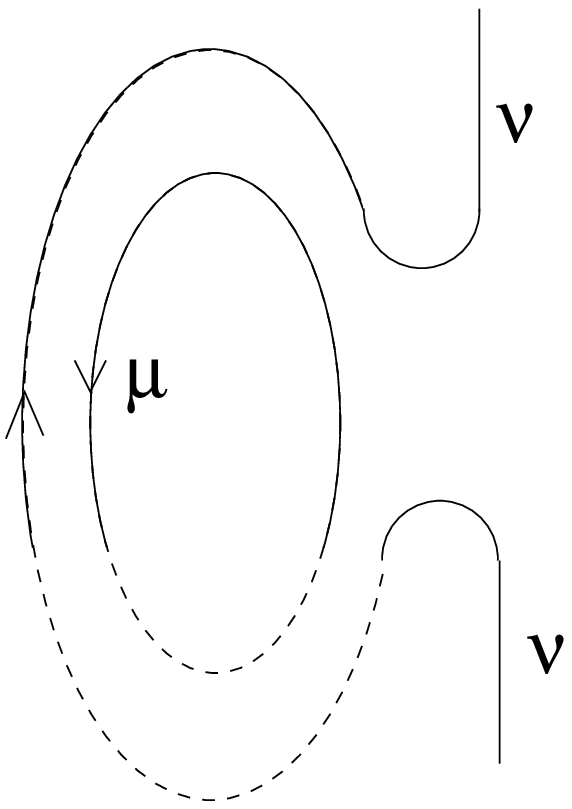}\, 
\end{eqnarray*}
In the first and third  equalities we use the fusion formula,
the second equality holds  by  isotopy.

\end{proof}

\noindent
  A more general statement is shown in \cite{AB}.

\begin{lem}(Killing property)
Suppose that $\la\omega\ra$ is nonzero.
 Let $\l\in {\Gamma (A)}$, then  
 the following morphism 
is nonzero  in $A$ if and only if $\l$ is transparent. 
$$\psdiag{7}{21}{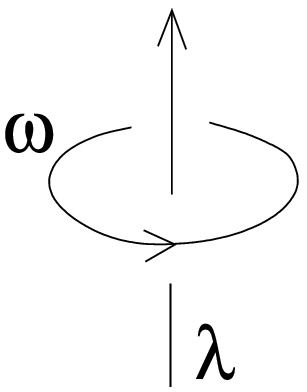}$$
\end{lem}

\begin{proof}
 If $\l$ is transparent, then this morphism
is equal to $\la\omega\ra\one\!_\l$, which is nonzero.
Conversely, if this morphism  is nonzero, it is equal to
$c \one\!_\l$ for some $0\neq c\in \fk$. Then,
for any $\n\in {\Gamma (A)}$, we have
$$\psdiag{8}{24}{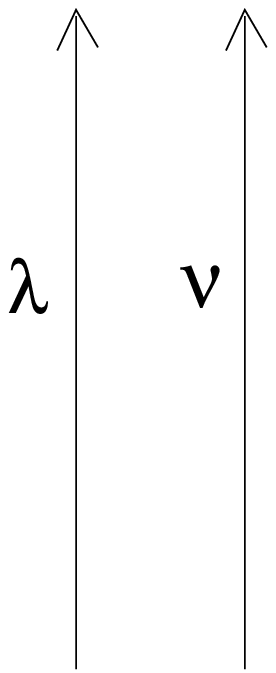}=c^{-1}\psdiag{8}{24}{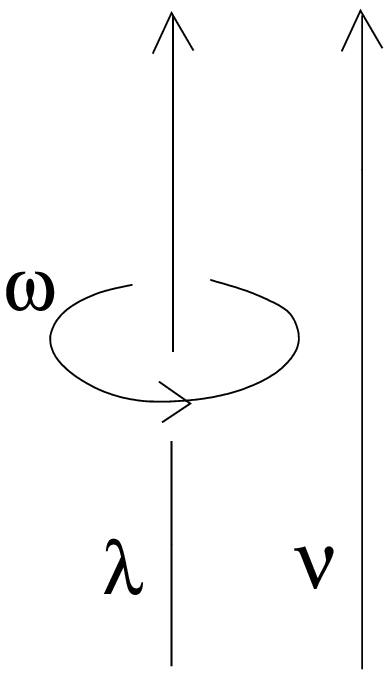}=
c^{-1}\psdiag{8}{24}{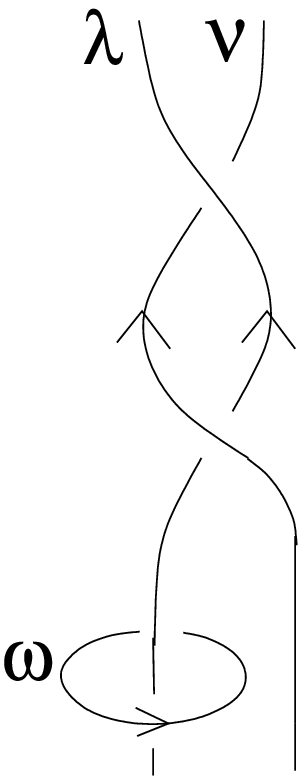}=\psdiag{8}{24}{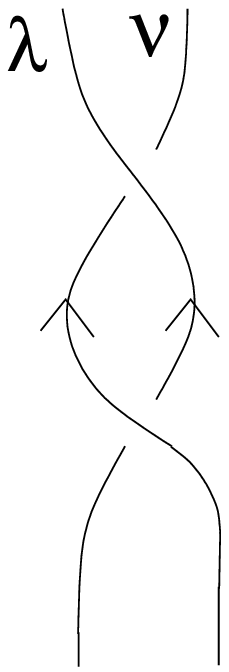}.$$

\noindent The second  equality holds by the sliding lemma.
\end{proof}

\begin{pro}\label{trmod}
A pre-modular category $A$ with $\la \omega \ra\neq 0$
which has no non-trivial transparent simple object is modular.
\end{pro}

\begin{proof}
We have to check the non-degeneracy axiom.
Let us denote by $\bar{S}$ the matrix whose $(\l,\m)$ entry is equal
to
 the value of the $0$-framed Hopf link with linking -1 
and coloring of the components $\l,\m$.
Then we have that
$$ \psdiag{7}{21}{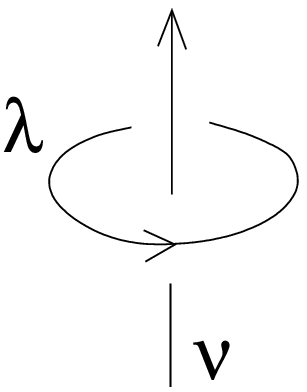}=\frac{S_{\l\n}}{\la \n\ra}
 \psdiag{7}{21}{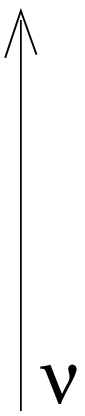}\ \ \text{ and } \ \ 
\psdiag{7}{21}{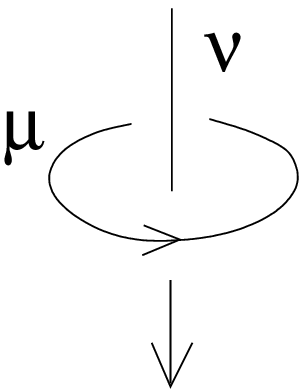}=\frac{\bar S_{\n\m}}{\la \n\ra}\psdiag{7}{21}{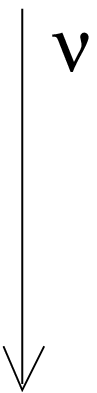}
\ .$$
 We deduce that the $(\l,\m)$ entry of the matrix $S\bar S$
is equal to the  invariant of the colored link depicted below.
$$\psdiag{8}{24}{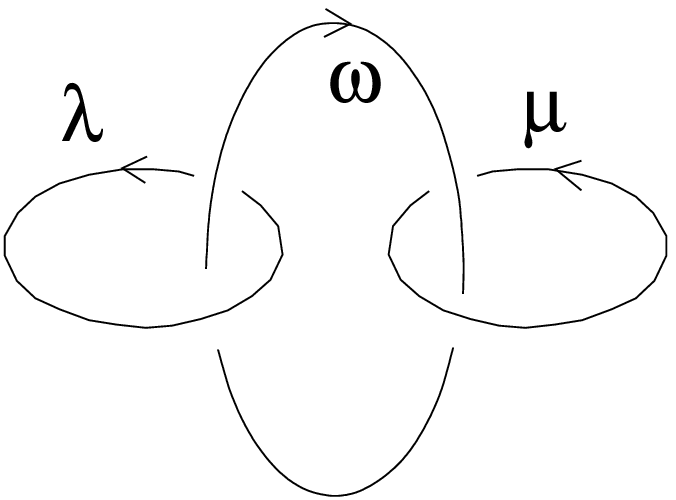}$$
 By using (\ref{tprod}) and the killing property we obtain
the formula 
 $$S\bar{S}=\la \omega\ra I\ ,$$
where $I$ is the identity matrix, which proves the invertibility of the
$S$ matrix.
\end{proof}

\subsection{Brugui\`eres' criterion}
A  process of constructing modular categories from   pre-modular ones
is called a  modularization.
Our reference for such construction is Brugui\`eres' work \cite{Bru}.
 See also \cite{Mu} for an analogous development in the context
of $*$-categories.
 Brugui\`eres considers abelian ribbon linear categories.
 Direct sums may be defined in a formal way, and a pre-modular category
with direct sums is an abelian category.
 From now on our pre-modular categories are supposed to be equipped
with direct sums (we add them if necessary) and hence are abelian.
\begin{dfn} 
A {\it modularization} of a pre-modular category $A$ is a 
modular category $\widetilde A$ together with a ribbon $\bf k$-linear functor
$F:A\to \widetilde A$ which is dominant,
i.e. any object of $\widetilde A$  
is a direct factor of $F(\l)$ for some $\l\in Ob(A)$.
\end{dfn}

\begin{dfn} A simple  object $\l$
 of  a pre-modular category $A$
is {\it bad} if for any $\m$ in a representative set of simple
objects ${\Gamma (A)}$, one has
$S_{\l\m}=\la \l\ra \la \m\ra $.
\end{dfn}

\begin{dfn}
For any $\l\in \Gamma(A)$, its twist coefficient 
$t_\l$ is defined by the  equality given below.
$$\psdiag{5}{15}{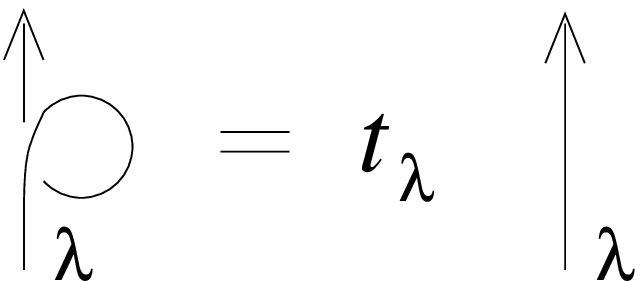}$$
\end{dfn}

\noindent
The following fact was claimed in Corollary 3.5 of \cite{Bru}.

\begin{thm}[Brugui\`eres' criterion]
Let  $\fk$ be  an algebraically closed field  of zero  characteristic.
 Then  an abelian pre-modular category $A$ over $\fk$
is modularizable if and only if any 
 bad object $X$
is transparent,  has twist coefficient $t_X=1$ and  
quantum dimension $\la X
\ra \in {\mathbb N}$.

If $A$ is modularizable, then its  modularization is unique up to
equivalence.
\end{thm}
\noindent
{\bf Remark.}
Clearly, any transparent object is bad.
If $\la\omega\ra\neq 0$, then any bad object  is transparent.
This follows from the killing property.
 Using this fact, Brugui\`eres statement can be slightly simplified  \cite{AB}.

\subsection{Modularization functor}

 We want now to describe the modularization functors explicitly.
 The main idea 
consists of  adding 
 morphisms to the pre-modular category, that  make  transparent simple
objects isomorphic  to the trivial one.

For the remainder  of this  section 
we consider a  pre-modular category $A$ with $\la\omega\ra\neq 0$,
whose
transparent simple objects 
have twist coefficient and quantum dimension equal to one.
 This corresponds to Brugui\`eres' particular case
\cite[Section 4]{Bru} and to M\"uger abelian case \cite[Section 5]{Mu}.
 The tensor product of two transparent simple objects is
then a transparent simple object, and  isomorphisms classes
of transparent simple object form a group $G$ under tensor multiplication.
 We will follow the description of the modularization functor
given in the proof of \cite[Lemma 4.3]{Bru}.
 As before, 
let  ${\Gamma (A)}$ be the representative set of simple objects of $A$.


 If $A$ is {\it self-dual} (i.e. any object is isomorphic to its dual),
then  
$G$ is isomorphic to
$({\mathbb Z}/2{\mathbb Z})^{|\cal T|}$, where $\cal T$ 
is the set of independent generators of $G$.
 This covers all cases considered in the next sections.

In general,  $G$ is isomorphic to
$\oplus^p_{i=1} {\mathbb Z}/k_i {\mathbb Z}$, $k_{i+1}|k_i$,
 and admits the following
presentation by generators and relations: 
$G\approx\{t_1,...,t_p;t_i^{k_i}=1, i=1,...,p\}$. 
 We fix, for each $i$,  a transparent simple object
representing the $i$th generator of $G$ 
and  denote it  by the same letter  $t_i$.  
Let  ${\cal T}=\{t_1,...t_p\}$ be  the   set of generating
transparent simple objects.
We denote by $G_{\mathcal{T}}$ the set of representatives
of $G$ defined by  $\mathcal{T}$, i.e.
$$G_{\mathcal{T}}=\{\otimes_it_i^{n_i}; t_i\in {\cal T}, 0\leq n_i<k_i\}\, .$$
Furthermore, we choose for each $i$
an  isomorphism $\Phi_i : t_i^{k_i}\approx {\rm trivial\;\; object}$.

Let us define a category $A'$ as follows. 
We set  $Ob(A')=Ob({A})$, we will however use the notation $F$
for the functor from $A$ to $A'$, 
and 
$$Hom_{A'}(F(X),F(Y)):=\oplus_{W\in G_\mathcal{T}}Hom_A(X,Y\otimes W)\ .$$
 For composition, we proceed as follows.
 Let $f\in Hom_A(X,Y\otimes W)$, $g\in Hom_A(Y,Z\otimes W')$ with 
$W, W'\in G_{\cal T}$.
 Since the objects of $G_{\cal T}$ are transparent, we get
a canonical isomorphism
$\Xi: Z\otimes W'\otimes W\to Z\otimes (\otimes_i t_i^{n_i})$.
We define
$F(g)F(f) :=\Xi (g\otimes \one\!_W)f$,
 if   $n_i<k_i$ for every $i$; 
otherwise  we compose the right hand side of the previous formula with
the isomorphisms $\one\!_{n_i-k_i}\otimes \Phi_i$ in order to reduce
the exponents.
 Associativity results from the property
\be \label{f2}\Phi_i\otimes \one\!_{t_i}=\one\!_{t_i}\otimes\Phi_i\ee
which is a consequence of 
\be \label{f1} \one\!_{t_i}\otimes \one\!_{t_i}=\psdiag{5}{15}{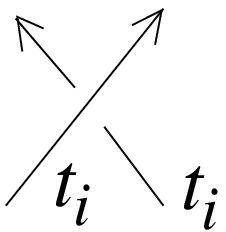} \ee
These are properties ($\mathcal{F}$) in \cite{Bru};
here we use that the transparent simple objects
are invertible,
 so that $t_i\otimes t_i$ is simple, and that their quantum dimensions
and twist coefficients are equal to one.

We define the category $\widetilde A$ as the idempotent completion of $A'$.
 It results from \cite[Section 4]{Bru} that $\widetilde A$ is a modularization
 of $A$. 

\v8
\noindent
{\bf Remark.} The category $\widetilde A$ is called sometimes a 
modular extension of $A$ by $G$. Analogously, a   modular extension
of $A$ by any subgroup $G'$ of $G$ can be constructed. This gives a 
pre-modular category  whose  group of transparent objects is $G/G'$.

\v8

The next problem is to construct a representative
 set ${\Gamma (\widetilde A)}$ of simple 
objects of $\widetilde A$.
 There is an action of the group $G$ on the set 
${\Gamma (A)}$ of simple objects of $A$ by tensor multiplication.  
For $X\in {\Gamma (A)}$, the dimension of 
$End_{\tilde A}(F(X))$ is equal to the order 
$d$ of the stabilizer subgroup $Stab(X):=\{g\in G; g\otimes X=X\}$.

If $Stab(X)$ is cyclic, then the algebra $End_{\widetilde A}(F(X))$
is abelian; it is isomorphic to the group algebra of $Stab(X)$,
and $F(X)$ decomposes in the category $\widetilde A$ into $d$ 
non-isomorphic simple objects.

In the non-cyclic case it can be shown (cf. \cite[Section 5]{Mu})
that $End_{\widetilde A}(F(X))$ is a twisted group algebra.
 The computation of the cocycle describing this twisted group
algebra has to be done. 

\subsection{Generalized ribbon graphs}
By Turaev's theorem \cite[Ch. I, Theorem 2.5]{Tu} 
the morphisms of a ribbon category $A$ can be represented by $A$-colored
ribbon graphs
with coupons. More precisely, there exists a functor from the
category $Rib_A$ of colored ribbon graphs to the category $A$ which respects 
the structures. We can extend the category $Rib_A$ by allowing tangles such 
that one of the ends of a band  colored with an object $t$ 
of  $\mathcal{T}$ is free. This means,
 it is connected neither to a coupon, nor to the source,
nor to the target. An example of such a tangle is depicted below.
It is considered as a morphisms from $Y$ to $X$.
$$\psdiag{8}{24}{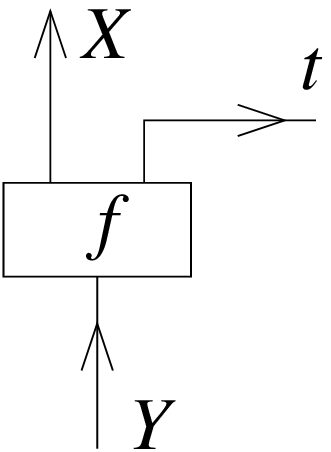}$$
This defines  the extended category $\widetilde{Rib}_A^\mathcal{T}$, which
is also a ribbon category. 
 We extend the invariant of closed colored graphs,
i.e the map $End_{Rib_A}(\text{trivial})\rightarrow \fk$ given by Turaev's functor,
in the following way.
 An extended closed colored graph is sent to zero, 
if the number of its free ends colored by $t_i$ is not divisible by $k_i$
for some $i$. Otherwise, it is sent to the invariant
of $Rib_A$ for a graph obtained by closing the free ends with $\Phi_i$.

Using the properties (\ref{f1}), (\ref{f2}), we can show that 
Turaev's functor
extends to a functor from $\widetilde{Rib}_A^\mathcal{T}$
to the modular category $\widetilde A$ which coincides
with the invariant above for closed morphisms.

\noindent{\bf Remark.} The modularization can be obtained from 
the ($\fk$-linear) category $\widetilde{Rib}_A^\mathcal{T}$
by first quotienting by negligible morphisms
(using the invariant $End_{Rib_A}(\text{trivial})\rightarrow \fk$
described above)
and then completing with idempotents.
 Direct sums are not needed here. This process was sketched in
\cite{Bl}.

\section{Idempotents of BMW algebras}

\subsection{Kauffman skein relations.}

Let  $M$ be a 3-manifold (possibly with a given finite set $l$ of  
points on the boundary, and a nonzero tangent 
vector at each point).
Let  $\fk$ be  a field
containing  the nonzero elements $\a$ and  $s$ with $s^2\neq 1$. 

 We denote by $\S(M)$  (resp. $\S(M,l)$)
the $\fk$-vector space freely generated by
links in $M$ (and 
tangles  in $M$ that meet $\partial M$ in $l$)
modulo 
the Kauffman  skein relations:
$$
\psdiag{3}{9}{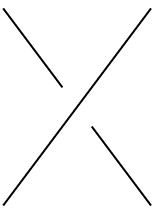}\;-\;\psdiag{3}{9}{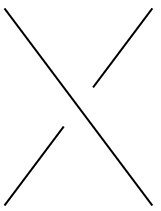} = \;(s-s^{-1})\;
\left(\;\,\psdiag{3}{9}{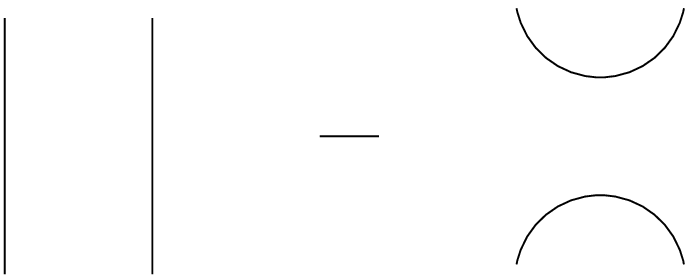}\;\,\right) 
$$
$$
\psdiag{3}{9}{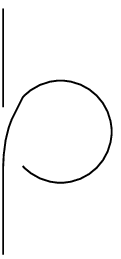}\;=\;\a\;\; \psdiag{3}{9}{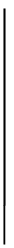}\;,\;\;\;\;\;
\psdiag{3}{9}{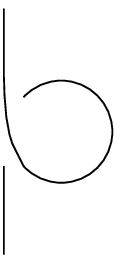}\;=\;\a^{-1} \;\;\psdiag{3}{9}{1}$$

$$L\;\amalg \;\bigcirc\; = \;\left(\frac{\a-\a^{-1}}{s-s^{-1}}+1
\right)\;\, L  .$$

\noindent
We call $\S(M)$ the {\it skein module}  of $M$. 
For example, $\S(S^3)\cong {\fk }$.
 
\subsection{Birman-Murakami-Wenzl  category}

The  Birman-Murakami-Wenzl (BMW)
 category $K$ is defined as follows. An object of $K$ 
is a standard oriented disc $D^2\subset \mathbb{C}$ equipped with a finite set of points and a nonzero tangent 
vector at each point. Unless otherwise specified, we will  use the second 
vector of the standard basis (the vector $\sqrt{-1}$ in complex notation).
 If $\beta=(D^2,l_0)$ and $\gamma=
(D^2,l_1)$ are two such objects,
the module $Hom_K(\beta,\gamma)$ is defined as the  skein module
 $\S(D^2\times [0,1],l_0\times 0\amalg
l_1\times 1)$.  Composition is given by stacking of  cylinders.
 We will use the notation 
$K(\beta,\gamma)$ for $Hom_K(\beta,\gamma)$
and $K_\beta$ for $End_K (\beta)$.
 The tensor product is defined by using 
$j=j_{-1}\amalg j_1:{ D}^2\amalg { D}^2\hookrightarrow { D}^2$,
 where, for $\epsilon=\pm 1$,  $j_\epsilon :{ D}^2\hookrightarrow { D}^2$
is the embedding which sends $z$ to $\frac{\epsilon}{2}+\frac{1}{4}z$.

The BMW category is a $\fk$-linear
ribbon  category. 
As before, we denote by  $\la f\ra \in \fk$  the quantum
trace  of $f\in K_\beta$.
The  BMW categories defined using the  parameters
$(\a,s)$
and $(\a,-s^{-1})$ are isomorphic.

Let us denote
by  $n$  the object of $K$
formed with  the $n$ points $\{ (2j-1)/n\; -1  ;
j=1,...,n \}$  equipped with
the  standard vector. 
 Composition in the category $K$ provides a $\fk$-algebra structure
on  $K_n=End_K(n)$, and we get  the Birman-Murakami-Wenzl (BMW)  algebra.

The BMW algebra $K_n$
 is a deformation of the Brauer algebra (i.e. the centralizer
algebra of the semi-simple  Lie algebras of type B,C and D).
 It  is known to be generically semi-simple 
and its simple components correspond to the partitions $\l=(\l_1,
...,\l_p)$ with $|\l|=\sum_i \l_i=n-2r$, $r=0,1,...,[n/2]$. 

\subsection{Idempotents} Let $\l$ be a partition.
 We denote by $\Box_\l$ the object of $K$ formed with one 
point for each cell  of the Young diagram associated with $\l$.
 If $c$ has index $(i,j)$
($i$-th row, and $j$-th column),
then the corresponding point in ${D}^2$ is
$\frac{j+i\sqrt{-1}}{n+1}$.
 In \cite{BB} we have constructed  minimal 
idempotents $\tilde y_\l \in K_{\Box_\l}$. 
Let us recall their main properties in the generic case
(i.e. with $\fk=\mathbb{Q}(\alpha,s)$).

\v8
\noindent{\em Branching formula:} 
\be \tilde y_\lambda\otimes \one\!{}_1=
\sum_{\genfrac{}{}{0pt}{2}{\lambda\subset \mu}{|\mu|=|\lambda|+1}}
\psdiag{11}{33}{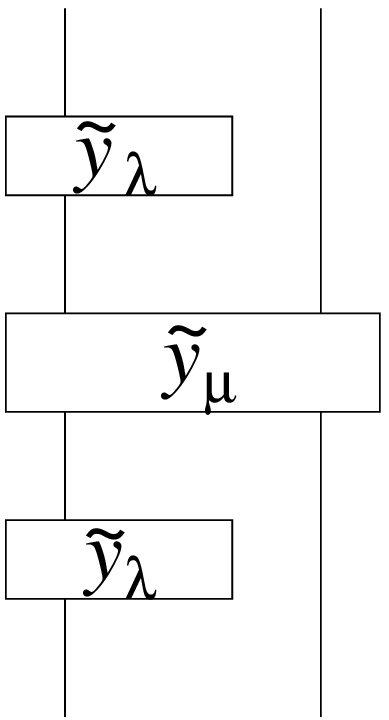}
+
\sum_{\genfrac{}{}{0pt}{2}{\mu\subset \lambda}{|\mu|=|\lambda|-1}}
\frac{\la \m\ra}{\la\l\ra}\;\;\; \psdiag{11}{33}{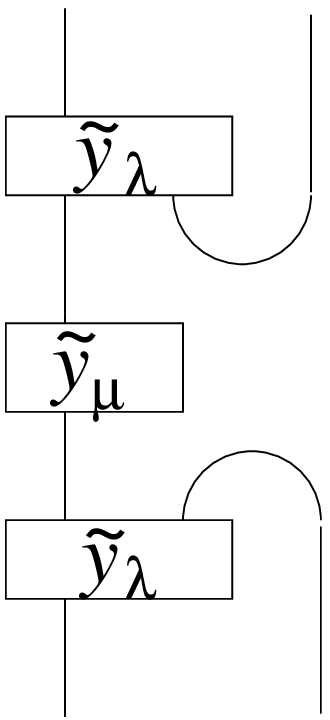} \ee

\noindent
 Here standard isomorphisms are used, in the first tangle
between $\Box_\l\otimes 1$ and $\Box_\m$, in the second tangle
 between $\Box_\m\otimes 1$ and $\Box_\l$.
The second tangle times $\frac{\la \m\ra}{\la\l\ra}$
will be further denoted by $\tilde y_{(\l,\m)}$.
Note that  the quantum dimension  $\la\l\ra$ is nonzero in the generic case. 

\v8
\noindent
{\em Braiding coefficient:} 
 Let $i)\;\m- \l=c$  or $ii)\; \l-\m =c$,  
where the cell $c$ has coordinates $(i,j)$. Let $cn(c)$ be the content
of the cell $c$: $cn(c)=j-i$. Then

\be\label{bc}
i)\;\;\;\;
\psdiag{15}{45}{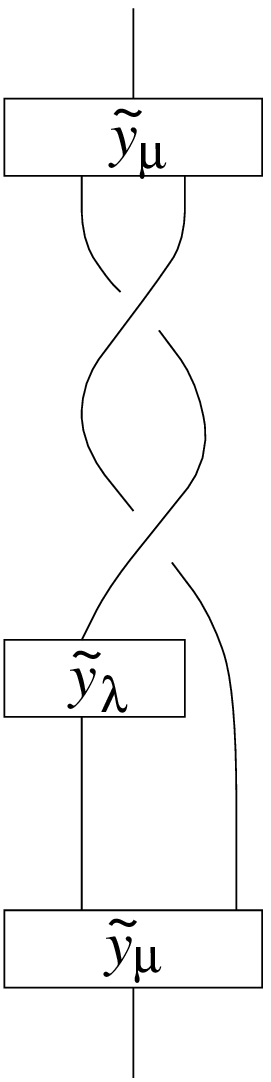}= s^{2cn(c)} \tilde{y}_\m ;\;\;\;\;\;
ii)\;\;\;\; 
\psdiag{15}{45}{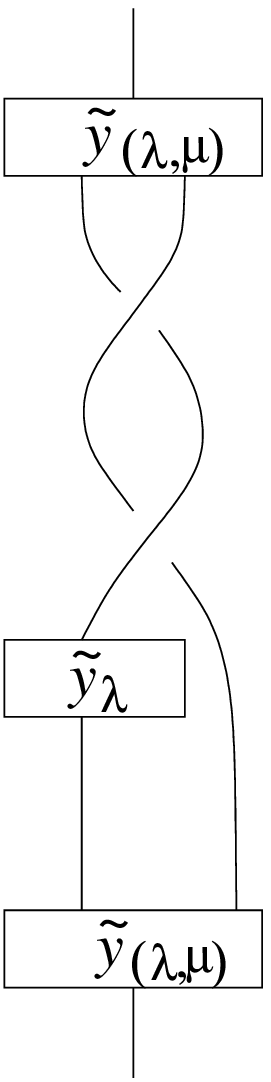}= \alpha^{-2} s^{-2cn(c)} \tilde{y}_{(\l,\m)} .\ee

\v8
\noindent
{\em Twist coefficient:}
 A positive $2\pi$-twist 
of $|\l|$ lines  with  $\tilde y_\l$ inserted
contributes  the factor 
$\a^{|\l|}s^{2\sum_{c\in \l} cn(c)}$. 
\be \psdiag{6}{18}{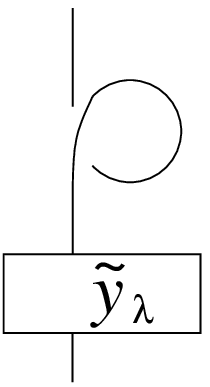}=\a^{|\l|}s^{2\sum_{c\in \l} cn(c)}
\psdiag{6}{18}{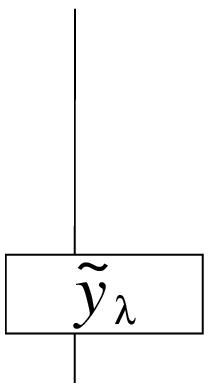}\ee

\v8
\noindent
{\em Quantum dimensions:}
Let $n\in {\mathbb Z}$, we set 
$$[n]_\alpha=\frac{\a
s^n-\a^{-1}s^{-n}}{s-s^{-1}}, \;\;\;\; [n]=\frac{s^n-s^{-n}}{s-s^{-1}}.$$
Then the quantum dimension of $\l$ is given by the following formula
\be \label{wen}\la \lambda\ra={\la \lambda\ra}_{\a,s}= 
\prod_{(j,j)\in\l}\frac{[\lambda_j-\l^\vee_j]_\alpha 
+[hl(j,j)]}{[hl(j,j)]} 
\underset{i\neq j}{\prod_{(i,j)\in\l}}\frac{[d_\l(i,j)]_\alpha}{[hl(i,j)]}.\ee 
Here, $hl(i,j)$ denotes  the hook-length of the cell 
$(i,j)$, i.e. $hl(i,j)=\l_i+\l^\vee_j-i-j+1$, $\l^\vee_i$ is the length
of the $i$-th column of $\l$ 
and $d_\l(i,j)$ is defined by 
$$d_\l(i,j)=\left\{\begin{array}{lcl} 
\lambda_i+\lambda_j-i-j+1&\text{ if }&i\leq j\\ 
-\l^\vee_i-\l^\vee_j+i+j-1&\text{ if }&i>j\ . 
\end{array}\right.$$ 
Observe that
\be\label{sym}
{\la\l\ra}_{\a,s}={\la\l\ra}_{-\a,-s}={\la\l\ra}_{\a^{-1},s^{-1}}=
\la \l^\vee\ra_{\a,-s^{-1}}\;\; .\ee
The formula (\ref{wen}) was first proved by
Wenzl \cite[Theorem 5.5]{Wbcd}. 
If we define $d'_\l(i,j)$ by
$$d'_\l(i,j)=\left\{\begin{array}{lcl} 
\lambda_i+\lambda_j-i-j+1&\text{ if }&i<j\\ 
-\l^{\!\vee}_i-\l^{\!\vee}_j+i+j-1&\text{ if }&i\geq j , 
\end{array}\right.  $$
then we can write Wenzl's formula  as follows.
\be \label{wenzltwo}\la \lambda\ra= 
{\prod_{(i,j)\in\l}}\frac{\alpha^\frac{1}{2}s^{\frac{1}{2}d_\l(i,j)}
-\alpha^{-\frac{1}{2}}s^{-\frac{1}{2}d_\l(i,j)}}
{s^{\frac{1}{2}hl(i,j)}-s^{-\frac{1}{2}hl(i,j)}}
{\prod_{(i,j)\in\l}}\frac{\alpha^\frac{1}{2}s^{\frac{1}{2}d'_\l(i,j)}
+\alpha^{-\frac{1}{2}}s^{-\frac{1}{2}d'_\l(i,j)}}
{s^{\frac{1}{2}hl(i,j)}+s^{-\frac{1}{2}hl(i,j)}}\ee  


\subsection{Idempotents in the non-generic case}\label{nongen}
By a non-generic case we understand a choice of parameters in 
the field $\fk$ such that $s$ is a root of unity,
or $\pm\a$ is a power of $s$.
A typical example is given by  roots of unity in a cyclotomic field.
 As in the generic case,
the idempotents $\tilde y_\l$ are obtained recursively
by lifting to the BMW category  the corresponding
idempotent $y_\l$ in the Hecke category.
 The minimal idempotent  $y_\l$
can be defined provided the quantum integers
$[m]$ are not zero for $m<\l_1+\l^\vee_1$, and $\tilde y_\l$
can further be obtained provided for {\em some}
$\mu\subset\l$, $|\mu|=|\l|-1$,
$\tilde y_\m$ is defined   and its quantum dimension is not zero.
 Under the above conditions, Wenzl path idempotent  \cite{Wbcd}
corresponding to a standard tableau $t$ with shapes
$\l(t)=\l$ and $\l(t')=\m$ is defined
and could be used here.
 The minimality property of the idempotent
$\tilde y_\l$ is 
$$\tilde y_\l K_{\Box_\l}\tilde y_\l= \fk \tilde y_\l\ .$$
The generic formulas of the previous subsection hold
provided they make sense.
 In particular
the branching  formula is valid provided the 
minimal idempotents exist for all diagrams obtained from $\l$ by 
adding one cell.

We will consider in the following the case where 
$\pm\a$ is a power of $s$, and discuss which idempotents are
obtained depending if 
$s$ is a root of unity or not. As explained in \cite{Wbcd},
in this case if we quotient out  the BMW algebra by
negligible morphisms (the annihilator of the trace),
 then we get a semi-simple algebra.

 If neither $\a$, nor $-\a$  are   powers of $s$
but $s$ is a root of unity, then  we obtain minimal idempotents corresponding to partitions
$\l$ with $\l_1+\l^\vee_1<l+1$, where $l$ is the order of $s^2$.
 These diagrams are called $l$-regular in \cite{Wbcd}.
If we consider a diagram $\mu$ with $\m_1+\m^\vee_1=l+1$,
obtained from an $l$-regular diagram by adding one cell,
then the generic element $\tilde Y_\mu=[l]\tilde y_\mu$ 
still can be defined and has  nonzero trace.
This element  satisfies
$\tilde Y_\mu K_{\Box_\mu}\tilde Y_\mu=0$, since $[l]=0$
in our specialization.
\begin{lem}
The element $\tilde Y_\mu$
belongs to the radical of the algebra $ K_{\Box_\mu}$
(the intersection of the maximal left ideals).
\end{lem}
\begin{proof}
Let $J$ be a maximal left ideal of  $ K_{\Box_\mu}$.
Suppose that $J$ does not contains
$\tilde Y_\mu$, then, using maximality of $J$, we get
that the left ideal $J+K_{\Box_\mu}\tilde Y_\mu$
is equal to $K_{\Box_\mu}$.
 We further 
have that $\one\!_{\Box_\mu}=j+a\tilde Y_\mu$, $j\in J$,
 $a\in K_{\Box_\mu}$,
 and so  $\tilde Y_\mu=\tilde Y_\mu j+\tilde Y_\mu a\tilde Y_\mu=\tilde Y_\mu j$
is in the ideal $J$, which contradicts the hypothesis.
\end{proof} 

This shows that the algebra $K_{\Box_\mu}$ is not semi-simple in this case
and 
if we quotient out   by
negligible morphisms we will still have a non semi-simple algebra.


\section{ The completed BMW categories }\label{trois}
In this section we define the completed BMW category and discuss
 specializations of parameters for which
the quotient of the completed BMW category
by negligible morphisms is a pre-modular category.

\subsection{Completed BMW categories}

 Let $\mathcal{C}$ be a set of Young diagrams,
 such that the corresponding minimal idempotents  exist.  
 This means that for each element of $\mathcal C$ the conditions described in 
Section \ref{nongen} are satisfied.
 In each case considered further this set will be the maximal set
in which the recursive construction of the idempotents $\tilde y_\l$
works (this set corresponds to the {\em affine Weyl alcove}
in the quantum group description).

 We define the  {\it completed BMW category} $K^\mathcal{C}$ as follows.
An object of $K^\mathcal{C}$
is an  oriented disc $D^2$ equipped with a finite set
 of  points, with a trivialization of the tangent space at each point
(usually the standard one), labeled with diagrams 
from $\mathcal{C}$.   
Let $\beta=(D^2,l)=(D^2;\l^{(1)},...,\l^{(m)})$ be such an object.
Then its expansion $E(\beta)=(D^2,E(l))$
is obtained by embedding the object $\Box_{\l^{(i)}}$
in a neighborhood of the  point labeled by  $\l^{(i)}$,
according to the trivialization.
The tensor product $\tilde y_{\l^{(1)}}\otimes...\otimes\tilde y_{\l^{(m)}}$
defines an idempotent $\pi_\beta\in K_\beta$. We define 
$
Hom_{K^\mathcal{C}}(\beta, \gamma):= \pi_\beta K(E(\beta), E(\gamma))\pi_\gamma\, .$
We will use the notation $K^{\mathcal{C}}(\beta, \gamma)$
and $K^{\mathcal{C}}_\beta$ similarly as in $K$.

The duality extends to $K^\mathcal{C}$, and we obtain
 again a $\fk$-linear ribbon category.
 Observe that the dual of an object 
is isomorphic to itself in a non-canonical way.

The equality of the categories $K$
 for the parameters $(\a,s)$ and $(\a,-s^{-1})$
extends to an isomorphism between the categories
$K^\mathcal{C}$ and $K^{\mathcal{C}^\vee}$, where
$\mathcal{C}^\vee$ is obtained from $\mathcal{C}$
by transposition of diagrams (i.e. exchange of rows and columns).
 For further discussion of duality, it is useful
to note that this change of the parameter $s$
switches a primitive $l$th root of unity, 
into a primitive $2l$th root of unity if $l$ is odd.

We denote by $\l$ the object of $K^\mathcal{C}$
formed by a disc with the origin labeled by $\l$.
 The minimality property of the idempotent
$\tilde y_\l$ implies that $\l$ is a simple object
in $K^\mathcal{C}$.

Recall that  a morphism  $f\in K^{\mathcal{C}}(\a,\beta)$ 
 is  negligible if
for any $g\in K^{\mathcal{C}}(\beta,\a)$ one has $\la fg\ra=0$.
 Negligible morphisms form a tensor ideal in the category, and
we obtain a quotient  $K^{\mathcal{C}}/Neg$ which is
 a $\fk$-linear ribbon category. The duality axiom is trivially 
satisfied here.
 Our  aim is to discuss in which case this quotient
category happen to be  pre-modular.

We first consider the { generic case}.
 Here the set $\mathcal{C}$
contains  all Young diagrams.  We see from the branching formula that
the completed category is semi-simple. Isomorphism classes
of simple objects correspond to all Young diagrams, so that
the category is not pre-modular. Moreover, from the braiding formula
(\ref{bc}) we see that there is no non-trivial
transparent simple object, so that we could not get a modularization
even if we would consider an extended version of Brugui\`eres' procedure.

We already have considered in Section \ref{nongen} the case where
 $s$ is a root of unity, but neither $\alpha$ nor $-\a$ is a power of $s$.
 Here the quotient of the idempotent 
completed category by negligible morphisms will not be semi-simple,
because some endomorphism algebras are not.

We will now consider  the  specializations where
$\pm \a$ is a power of $s$.
 Recall that $1^{\N+1}$ and $\K+1$ denotes the  column
and the  row Young diagrams with $\N+1$ and $\K+1$ cells,
respectively. 
Let us
 consider the following system of equations
 $\la 1^{\N+1}\ra=0$ 
 and $\la \K+1 \ra=0$, 
with $\N$ and $\K$  minimal. Note that, if $\pm \a$ is a power of $s$, then at least
one of these two equations has a solution.
 The first one is equivalent to 
$\a=-s^{2\N+1}$ or $\a=\pm s^{\N-1}$.  
 We have to consider 4 cases. 
\begin{verse} 
Case $C_n$: $\a=-s^{2n+1}$ ($\N=n$),\\ 
Case $B_n$: $\a= s^{2n}$ ($\N=2n+1$),\\ 
Case $B_{-n}$: $\a=- s^{2n}$ ($\N=2n+1$),\\ 
Case $D_n$: $\a= s^{2n-1}$ ($\N=2n$),\\ 
\end{verse}
 The interpretation of the notation $C_n$, $B_n$, $D_n$ is that the  
given specialization
of the Kauffman polynomial is obtained by using the fundamental 
representation of the 
corresponding quantum group. The specializations $B_n$ and $B_{-n}$ are 
similar, but
they are not equivalent; one should think of
the fundamental object in the $B_{-n}$ specialization as the deformation
 of the 
fundamental representation of $so(2n+1)$, with negative dimension $-2n+1$.

 The discussion of the  equation $\la \K+1\ra$ is similar. 
 Note that quantum dimensions are unchanged 
 if we replace $s$ by $-s^{-1}$ and interchange rows with columns. 
 Here are the four  cases. 
 \begin{verse} 
Case $C_k$: $\a=s^{-2k-1}$ ($\K=k$),\\ 
Case $B_k$: $\a= s^{-2k}$ ($\K=2k+1$),\\ 
Case $B_{-k}$: $\a=- s^{-2k}$ ($\K=2k+1$),\\ 
Case $D_k$: $\a=- s^{-2k+1}$ ($\K=2k$),\\ 
\end{verse} 
  We observe that, if $\la 1^{\N+1}\ra= 
\la \K+1\ra=0$ for some $\N$, $\K$, then $s$ is a root of unity.
 We will consider the four cases corresponding
to the vanishing of $\langle 1^{\N+1}\rangle$, and then,
according to the order of $s^2$,
combine them with  the condition corresponding
to the lowest $\K$ for which  $\la \K+1\ra$ vanishes.

The  cases $\a=\pm 1$, $\a=-s$ and $\a=s^{-1}$ will be excluded
from the general discussion given in the next subsections.
If $\a=\pm 1$ we get a category with two simple objects:
the trivial object and $\l=1$. The second object is transparent
and the category is modularizable iff $\a=1$.
The corresponding link invariant is trivial.
If $\a=-s$ or $\a=s^{-1}$, then the Kauffman polynomial is zero.

The case $\a=s$ (resp. $\a=-s^{-1}$) will be included in the general
discussion and give the categories $D^{1,k}$, $D\!B^{1,k}$
and $D\!B^{1,-k}$ (resp. $D^{k,1}$, $B\!D^{k,1}$
and $B\!D^{-k,1}$). Note that the corresponding invariant
of a link $L=(L_1,\dots,L_m)$ is equal to
$2^{\sharp L}s^{\sum_i L_i. L_i}$.
Here $\sharp L=m$ is the number of components,
and $L_i.L_i$ is the self linking number (the framing coefficient).
 The category is modularizable if $s$ is either a primitive root of order
$2l$, $l$ even, or a primitive root of odd order $l$.
 One can show that the corresponding invariants of $3$-manifolds are those
known as the $U(1)$ invariants \cite{MOO}.


\subsection{The symplectic case}
In this subsection  let  $\a=-s^{2n+1}$ , $n\geq 1$.
 (For $n=1$ the specialized Kauffman polynomial is
 the Kauffman bracket, and we will recover
 the TQFT's obtained in \cite{BHMV}.)

 If $s$ is generic, then we can construct the idempotent
$\tilde y_\l$ for $\l$ in the set
$$
\bar \Gamma(C_{n})
=\{\l;  {\l}^\vee_1\leq n+1,\  \l^\vee_2\leq n\}\ ,$$
and $\l$ has non-vanishing quantum dimension (see formula (\ref{wenzltwo}))
if it belongs to
$$\Gamma(C_{n})
=\{\l;  {\l}^\vee_1\leq n\}\ .$$
From the branching formula we get that the category 
$K^{{\Gamma}(C_{n})}/Neg$ is semi-simple; we will give more details
in the proof of Proposition \ref{premod}.
 A representative set of simple objects is
the infinite set $\Gamma(C_{n})$,
so that the category is not pre-modular.

The formula for the
quantum dimension can be simplified as follows  (see \cite[Prop. 7.6]{BB},
compare \cite{E}). 

\begin{pro}\label{qdimC}
 Let  $\a=-s^{2n+1}$, with $s$ generic. Then,
for a partition  $\lambda=(\l_1,...,\l_n)$, we have
$$\la\l\ra=(-1)^{|\l|}
\prod_{j=1}^n\frac{[2n+2+2\lambda_j-2j]}{[2n+2-2j]}
\prod_{1\leq i<j\leq n}\frac{[2n+2+\lambda_i-i+\lambda_j-j]
[\lambda_i-i-\lambda_j+j]}
{[2n+2-i-j]
[j-i]}.
$$
\end{pro}

\noindent

Let us suppose now that $\a=-s^{2n+1}$  with $s^2$ a primitive $l$th root 
of unity and $l\geq 2n+1$. One can check that the above formula
for quantum dimensions is still valid  provided
$l\geq 2n+1$. The condition $l\geq 2n+1$ ensures that $1^{n+1}$ 
is the smallest column with vanishing quantum dimension.
 Note that for $l=2n+1$ we have $\a=\pm 1$,
and for $l={2n+2}$, we have $\a=-s$. In the following we 
 discuss the equation $\la \K+1 \ra=0$ with
$\K$ minimal  according to $l\geq 2n+3$. 
\begin{itemize}
\item
 If $l\geq 2n+4$ is even, then  $\K=l/2-n-1=k$, and
$\a=-s^{2n+1}=s^{-2k-1}$.
This will be the $C_n$-$C_k$ specialization.
\item
If $l\geq 2n+3$ is odd and $s^l=-1$, then $\K=2k+1$, 
$\a=-s^{2n+1}=s^{-2k}$.
This will be the
$C_n$-$B_k$ specialization.
\item
If $l\geq 2n+3$ is odd and $s^l=1$, then $\K=l-2n=2k+1$, 
$\a=-s^{2n+1}= -s^{-2k}$.
 This will be the $C_n$-${B}_{-k}$ specialization.
\end{itemize}

The specializations  $C_n$-$B_k$ and $C_n$-${B}_{-k}$ 
 are similar because of the symmetry
$(\a,s)\leftrightarrow (-\a,-s)$ for quantum dimensions.
 Note however that the twist coefficient is not preserved under
 this symmetry,
so that the modularization problems will be distinct.
 We will show that the $C_n$-$C_k$ and $C_n$-$B_k$ specializations
 lead to modular categories.

\v8\noindent  
{\bf  ${\bf  C}^{n,k}$  category.} 
Let us consider  the  $C_n$-$C_k$
 specialization of parameters with $n,k\geq 1$, i.e.
$\a=-s^{2n+1}=s^{-2k-1}$ and $s$ is a primitive
$2l$th root of unity with $l=2n+2k+2$.
 We will use the following sets of Young diagrams:
$$
\bar \Gamma(C^{n,k})
=\{\l; \l_1\leq k+1, \l_2\leq k, {\l}^\vee_1\leq n+1, \l^\vee_2\leq n
\}\ ,$$
$$\Gamma(C^{n,k})
=\{\l; \l_1\leq k, {\l}^\vee_1\leq n\}\ .$$
 We can construct the minimal  idempotent 
 for each $\l\in \Gamma(C^{n,k})$, since
the quantum dimensions of these objects given by Proposition 
\ref{qdimC} do not vanish.
Let $\l\in \Gamma(C^{n,k})$. If $\mu$ is obtained from
$\l$ by adding one cell, then $\tilde y_\mu\in \bar\Gamma(C^{n,k})$
can be constructed.
 Moreover, if 
 $\mu$ is not in  $\Gamma(C^{n,k})$, then $\langle \mu\rangle$
vanishes, and so $\tilde y_\mu$ is negligible.

 The category $C^{n,k}$ is defined as the quotient  
of the category
  $K^{\Gamma(C^{n,k})}$ by negligible morphisms.

\v8\v8
\noindent
{\bf  ${\bf CB}^{n,k}$ and ${\bf C B}^{n,-k}$ categories.}
 In the  case of the $C_n$-$B_k$ (resp. $C_n$-${B}_{-k}$)
specialization with $n,k\geq 1$ we have
$\a=- s^{2n+1}= s^{-2k  }$ and  $s$ is  a primitive $2l$th root of unity
(resp. $\a=- s^{2n+1}= -s^{-2k}$ and  $s$ is 
 a primitive $l$th root of unity),
$l=2n+2k+1$.
We proceed as above with
$$
\bar \Gamma(C\!B^{n,k})=\bar \Gamma(C\!{ B}^{n,-k})
=\{\l; \l_1+\l_2 \leq 2k+2,
\l^\vee_1 \leq n+1,\l^\vee_2\leq n \},$$
$$ \Gamma(C\!B^{n,k})=\Gamma(C\! { B}^{n,-k})=\{\l; \l_1+\l_2 \leq 2k+1,
\l^\vee_1\leq n\},$$ 
$$C\!B^{n,k}=K^{\Gamma(C\!B^{n,k})}/{\sf Neg}\, , \;\;\;
C\!B^{n,-k}=K^{\Gamma(C\!{ B}^{n,-k})}/{\sf Neg}\ .$$

\begin{pro}\label{premod}
For $n,k\geq 1$, the categories $C^{n,k}$,  
   $C\!B^{n,k}$ and $C\!{ B}^{n,-k}$ with representative 
sets of simple objects 
$\Gamma(C^{n,k})$,
 $\Gamma(C\!B^{n,k})$ and
$\Gamma(C\!B^{n,-k})$, respectively,
are pre-modular.
\end{pro}
\begin{proof}
 We have to prove the dominating property.
 The proof is the same in all cases, so we will use
the notation $\bar\Gamma$, $\Gamma$ for $\bar\Gamma(A)$, $\Gamma(A)$
where $A$ is one of the categories mentioned in the claim.
 It is enough to show that the identity morphism
of the  object $n$ decomposes using the simple objects
in $\Gamma$. This is done by induction on
$n$. For the step from $n$ to  $n+1$,
we have to decompose $\one\!_\l\otimes \oneun$, with
$\l\in\Gamma$. The key point is that
any diagram obtained from $\l$ by adding one cell
is in $\bar \Gamma$. Hence we have that the
branching formula holds and gives the required decomposition,
because the idempotents indexed by partitions in
$\bar \Gamma\setminus \Gamma$ are negligible.
\end{proof}

\subsection {The odd orthogonal case}
  We first consider the $B_n$ specialization  $\a=s^{2n}$.
 If $s$ is generic, then we can construct the idempotent
$\tilde y_\l$ for $\l$ in the set
$$
\bar \Gamma(B_{n})
=\{\l;  {\l}^\vee_1+ \l^\vee_2\leq 2n+2\}\ ,$$
and $\l$ has non-vanishing quantum dimension (see formula (\ref{wenzltwo}))
if it belongs to
$$\Gamma(B_{n})
=\{\l;  {\l}^\vee_1+ \l^\vee_2\leq 2n+1\}\ .$$
As we did before, we get that the category 
$K^{ \Gamma(B_{n})}/Neg$ is semi-simple.
 A representative set of simple objects is
the infinite set $\Gamma(B_{n})$,
so that the category is not pre-modular.
 
We have the following specialized formula for the quantum dimensions
(see \cite[ Prop. 7.6]{BB}).
\begin{pro}\label{qdimB} Let $\a=s^{2n}$, 
 with $s$ generic.
 For a partition  $\lambda=(\l_1,...,\l_n)$, we have
$$\la\l\ra=\prod_{j=1}^n\frac{[n+\lambda_j-j+1/2]}{[n-j+1/2]}
\prod_{1\leq i< j\leq n}\frac{[2n+\lambda_i-i+\lambda_j-j+1]
[\lambda_i-i-\lambda_j+j]}
{[2n-i-j+1]
[j-i]}.$$
\end{pro}
\noindent
In this case, the object $1^{2n+1}$ plays a special role.
\begin{lem}
Suppose that $\a=s^{2n}$, 
and $s$ is generic.
 Then the object $1^{2n+1}$ is transparent and it is the unique nontrivial 
transparent object in $\Gamma(B_{n})$. Its quantum dimension 
and twist coefficient are equal to one.
\end{lem}
\begin{proof}
 An object $\l\in \Gamma(B_{n})$ is transparent if and only
if for any (non-negligible) $\mu$ in the branching formula for $\l$,
the braiding coefficient 
is equal to one. 
Indeed, if all braiding coefficients are equal to one,
by summing over $\m$ the left hand sides and right hand sides
of (\ref{bc}) and applying the branching formula we have
$$\psdiag{7}{21}{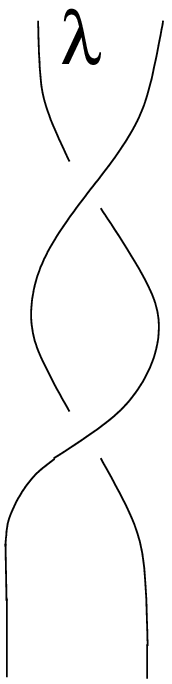}=\psdiag{7}{21}{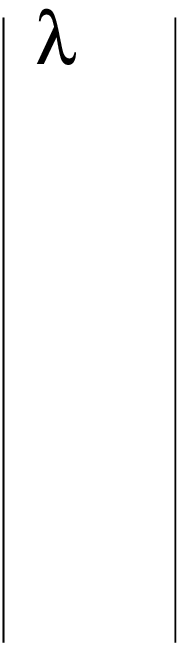}.$$  
Using this equality repeatedly we conclude that $\l$ is transparent.
Conversely, if $\l$ is transparent, its braiding coefficients are 
trivial.
 
The object
 $1^{2n+1}$ has only one braiding coefficient
 corresponding
to the removal of the last cell, and this coefficient is one.
(Two diagrams obtained by adding one cell to $1^{2n+1}$
are negligible.)
 It remains to check that any nontrivial $\l\in\Gamma(B_{n})$ distinct from
$1^{2n+1}$ has at least one braiding coefficient distinct from $1$.
 If $\mu$ is obtained from such $\l$
 by adding a cell in the first row, 
then $\la\mu\ra$ is not zero,
and the corresponding  braiding coefficient in formula
(\ref{bc}) is $s^{2\l_1}\neq 1$.
 For a column with $j$ cells,
 the generic quantum dimension reduces
to 
\be\label{col}
\la 1^j\ra=\frac{[0]_\a[-1]_\a...[2-j]_\a([1-j]_\a+[j])}{[j]!}.\ee
This gives for $1^{2n+1}$
 $$\la 1^{2n+1}\ra=\frac{[2n]\dots[1](0+[2n+1])}{[2n+1]!}=1\ .$$
The twist coefficient for $1^{2n+1}$ is $\a^{2n+1}s^{-2n(2n+1)}=1$.
\end{proof}
\begin{pro}\label{isogen}
In the category $K^{ \Gamma(B_{n})}/Neg$, \\
a) the object
$1^{2n+1}\otimes 1^{2n+1}$ is isomorphic to the trivial object;\\
b) 
the objects $1^{2n+1}\otimes \l$ and $\tilde \l$ are isomorphic,
where   $\l 
\in \Gamma(B_n)$, and  $\tilde \l$ is  the Young diagram such that
$\l_1^\vee+\tilde \l_1^\vee=2n+1$ and $\l_j^\vee=\tilde \l_j^\vee$ for $j>1$,
\end{pro}
\begin{proof}
In the semi-simple category $K^{ \Gamma(B_{n})}/Neg$
 we can decompose the identity of the object
 $1^{2n+1}\otimes 1^{2n+1}$ as we did in formula
(\ref{tprod}).
$$\sum_{\m}\sum_{i}\la\m\ra
\ \ \psdiag{11}{33}{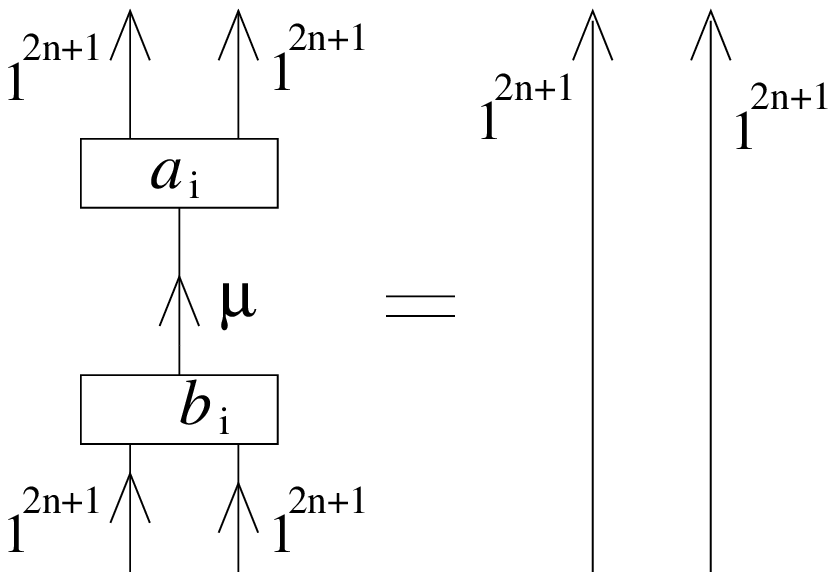}\ .$$
 Here all simple subobjects $ \m$ are transparent 
and hence have dimension $1$.
 By comparing the dimensions we see that there is only one such $\mu$ 
with multiplicity $1$.
 It should be trivial, because the duality gives a nonzero morphism
from the trivial to $1^{2n+1}\otimes 1^{2n+1}$.
 We deduce that this duality morphism is an
 isomorphism, which establishes $a)$.

We consider the   morphism from  $1^{2n+1}\otimes \l$ to $\tilde \l$
depicted below: the strings corresponding to the points
in (the expansion of)  $1^{2n+1}$ are joined to the first columns,
the points which are not in the first column of $\l$
and $\tilde \l$ are joined directly.
$$\psdiag{15}{40}{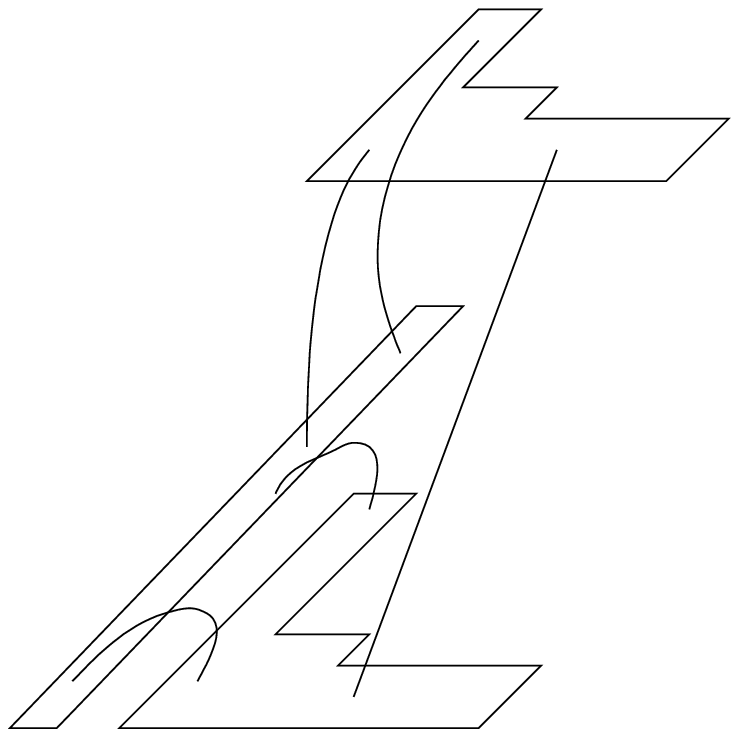}$$
One wants to show that this morphism is nonzero.
   We first consider the case where
 $\l=1^j$ has only one column.
Let $f\in Hom(1^j\otimes 1^{2n+1}, 1^{2n+1-j})$ be the  morphism
as above
and $g\in Hom(1^{2n+1-j},1^j\otimes 1^{2n+1} )$ be its mirror image with
respect to the target plane. Then
  $\la gf\ra=\la 1^{2n+1}\ra=1$.
  In the general case, if we insert conveniently
 the isomorphism considered in the particular case
 between and $\one\!\!_{ 1^{2n+1}\otimes \l}$ and  $\one\!_{\tilde \l}$
 we obtain our nontrivial morphism.
\end{proof}

We suppose now that $\a=s^{2n}$, with $s^2$  a 
primitive  $l$th root of unity,  $l\geq 2n+1$.
  In the following we 
 discuss the equation $\la \K+1 \ra=0$,
$\K$ minimal.
 If $s$ has order $2n+1$ and $s^l=1$, then $\a=s^{-1}$ 
and the Kauffman polynomial is trivial.
\begin{itemize}
\item
 If $l\geq 2n+2$ is even, then $\K=l-2n+1=2k+1$,
$\a=s^{2n}=-s^{-2k}$;
this will be the $B_n$-${ B}_{-k}$ specialization.

\item
If $l\geq 2n+1$ is odd and $s^l=-1$, then $\K=l+1-2n=2k$,
$\a=s^{2n}=-s^{-2k+1}$; this will be
$B_n$-$D_k$ specialization.

\item
If $l\geq 2n+3$ is odd and $s^l=1$, then $\K=\frac{l-1}{2}-n=k$, 
$\a=s^{2n}= s^{-2k-1}$ will be the $B_n$-$C_k$ specialization.
\end{itemize}

\v8\v8 
\noindent
  {\bf  ${\bf  {B}}^{n,-k}$ 
  category.} 
Here we consider the $B_n$-${ B}_{-k}$ specialization 
($\alpha=s^{2n}=-s^{-2k}$) 
with  $n,k\geq 1$,
 $s$ is a primitive $2l$th root of unity,
$l=2n+2k$. 
Let 
$$ \Gamma({B}^{n,-k})=
\{\l; \l_1+\l_2 \leq 2k+1,
\l^\vee_1 +\l^\vee_2\leq 2n+1\}
 \ .$$ 

We can define idempotents for
any $\l\in  \Gamma(B^{n,-k})$, and they have nonzero
quantum dimension. 
 Our general procedure give some more idempotents
whose dimension vanishes, namely for each 
$\l\in \bar \Gamma(B^{n,-k})\setminus \Gamma(B^{n,-k})$ with
$$\bar\Gamma({B}^{n,-k})
=\{\l; \l_1+\l_2 \leq 2k+2,
\l^\vee_1 +\l^\vee_2\leq 2n+2, \l_1+\l^\vee_1\leq 2n+2k\}$$ 
we have $\la \l\ra=0$.
We define the  category $B^{n,-k}$ 
as the quotient  of the category $K^{\bar \Gamma(B^{n,-k})}$
by negligible morphisms.
\begin{pro}
The category $B^{n,-k}$
is pre-modular.
\end{pro}
\begin{proof}
Let $\tilde \Gamma(B^{n,-k})
=\Gamma (B^{n,-k})\cup \{1^{2n+1}\otimes 2k+1\}$.
We show that 
$\tilde \Gamma(B^{n,-k})$
is a
set of dominating simple objects.
 As in the proof of Proposition \ref{premod}, 
we  decompose 
the tensor products $\one\!_W\otimes \one\!_1$,
for $W\in  \tilde \Gamma(B^{n,-k})$.
The sublte point here
is that some idempotent in the
 branching formula for the partition $L=(2k,1^{2n-1})\in \Gamma(B^{n,-k})$
(i.e. $L_1+L^\vee_1=2n+2k$) is  missing. 
 We  will  avoid this difficulty by
using the isomorphism in Proposition \ref{isogen}
which still holds  for 
$\l\in \Gamma(B^{n,-k})$. 

More precisely,
if $W=\l$ is in $\Gamma(B^{n,-k})\setminus \{L\}$, then
 the branching formula applies.
If $W=L$, then we use the isomorphism between $L$
and $1^{2n+1}\otimes 2k$ and we get a decomposition
of $L\otimes 1$ with subobjects
$(2k-1,1^{2n})$, $(2k,1^{2n-1})$ and $1^{2n+1}\otimes 2k+1$.
If $W=1^{2n+1}\otimes 2k+1$, then we get
an isomorphism between $1^{2n+1}\otimes 2k+1\otimes 1$ and
$L$.
\end{proof}


\v8
\noindent {\bf  ${\bf BD}^{n,k}$ category.}
 For the $B_n$-$D_k$ 
 specialization  with $n,k\geq 1$,
 we put  $l=2n+2k-1$, $s$ is a 
primitive root of unity of order $2l$, and $\a= s^{2n}=- s^{-2k+1  }$. Let
$$ \Gamma(B\!D^{n,k})= \{\l; \l_1+\l_2 \leq 2k,
\l^\vee_1+\l^\vee_2\leq 2n+1\} \ .$$
 We define the  category $B\!D^{n,k}$ 
 and prove pre-modularity as we did above.

\v8
\noindent {\bf  ${\bf BC}^{n,k}$ category.} 
The category $BC^{n,k}$ for $n,k\geq 1$  with parameters $(\a,s)$
is isomorphic to the category $C\!B^{k,n}$
with parameters $(\a, -s^{-1})$. The isomorphism sends
any simple object $\l$ to $\l^\vee$. 
 The representative set of simple objects is 
$\Gamma(BC^{n,k})=\{\l;\l^\vee\in \Gamma (C\!B^{k,n})\}$.

\v8\noindent
{\bf The specialization ${\bf B}_{-n}$.}
Let us consider  the case 
$\a=-s^{2n}$. If $s$ is generic,
we have $\Gamma(B_n)=\Gamma(B_{-n})$.
The object $1^{2n+1}$ remains transparent,
but its twist coefficient is $(-1)$. Therefore, the categories 
we get here will  be non-modularizable.

Let us suppose that $s^2$ is a primitive root of unity
of order $l\geq 2n+1$. Then we have to consider the following cases.
\begin{itemize}
\item
 If $l\geq 2n+2$ is even, then $\K=l-2n+1=2k+1$,
$\a=-s^{2n}=s^{-2k}$;
this will be the $B_{-n}$-${ B}_{k}$ specialization.

\item
If $l\geq 2n+1$ is odd and $s^l=1$, then $\K=l+1-2n=2k$,
$\a=-s^{2n}=-s^{-2k+1}$; this will be
$B_{-n}$-$D_k$ specialization.

\item
If $l\geq 2n+3$ is odd and $s^l=-1$, then $\K=\frac{l-1}{2}-n=k$, 
$\a=-s^{2n}= s^{-2k-1}$ will be the $B_{-n}$-$C_k$ specialization.
\end{itemize}

The categories $B^{-n,k}$, $BD^{-n,k}$ and $BC^{-n,k}$ with $n,k\geq 1$
 can be constructed 
analogously to the previous case.
We have $\Gamma(B^{-n,k})=\Gamma(B^{n,-k})$,
$\Gamma(BD^{-n,k})=\Gamma(BD^{n,k})$ and
$\Gamma(BC^{-n,k})=\Gamma(BC^{n,k})$.

\subsection{The even orthogonal case}
In this subsection  we suppose that $\a=s^{2n-1}$, $n\geq 1$.
 If $s$ is generic, then we can construct the idempotent
$\tilde y_\l$ for $\l$ in the set
$$
\bar \Gamma(D_{n})
=\{\l;  {\l}^\vee_1+ \l^\vee_2\leq 2n+1\}\ ,$$
and $\l$ has non-vanishing quantum dimension (see formula (\ref{wenzltwo}))
if it belongs to
$$\Gamma(B_{n})
=\{\l;  {\l}^\vee_1+ \l^\vee_2\leq 2n\}\ .$$
We get that the category 
$K^{ \Gamma(D_{n})}/Neg$ is semi-simple.
 A representative set of simple objects is
the infinite set $\Gamma(D_{n})$.

We have the following 
 specialized formula for the
quantum dimension. 

\begin{pro}\label{qdimD} Let $\a=s^{2n-1}$,  with $s$ generic.
 For a partition  $\lambda=(\l_1,...,\l_n)$, we have
$$\la\l\ra=
\prod_{1\leq i<j\leq n}\frac{[2n+\lambda_i-i+\lambda_j-j]
[\lambda_i-i-\lambda_j+j]}
{[2n-i-j]
[j-i]}\text{ if $\lambda_n=0$;}
$$
$$\la\l\ra=
2\prod_{1\leq i<j\leq n}\frac{[2n+\lambda_i-i+\lambda_j-j]
[\lambda_i-i-\lambda_j+j]}
{[2n-i-j]
[j-i]}\text{ if $\lambda_n\neq0$.}
$$
\end{pro}
\noindent
Suppose that $s^2$ is a primitive
$l$th root of unity with $l\geq 2n$.
We discuss the equation $\la \K+1 \ra=0$, $\K$ minimal.
\begin{itemize}
\item
 If $l\geq 2n$ is even, then $\K=l-2n+2=2k$,
$\a=s^{2n-1}=-s^{-2k+1}$;
this will be the $D_n$-$D_k$ specialization.
\item
If $l\geq 2n+1$ is odd and $s^l=1$, then $\K=l-2n+2=2k+1$,
$\a=s^{2n-1}=s^{-2k}$; this will be
$D_n$-$ B_k$ specialization.
\item
If $l\geq 2n+1$ is odd and $s^l=-1$, then $\K=l-2n+2=2k+1$, 
$\a=s^{2n-1}= -s^{2k}$ will be the $D_n$-$B_{-k}$ specialization.

\end{itemize}

 
\v8\noindent
{\bf  ${\bf D}^{n,k}$ category.} 
We  consider the $D_n$-$D_k$ specialization with  $n,k\geq 1$.
 Let
$$ \Gamma(D^{n,k})=\{\l; \l_1+\l_2 \leq 2k,
\l^\vee_1 +\l^\vee_2\leq 2n\}
\ .$$ 
We define the  category $D^{n,k}$ and prove pre-modularity
as above. The dominating set of simple objects is here
$\Gamma(D^{n,k})\cup \{1^{2n}\otimes 2k\}$.

\v8
\subsection{ The level-rank duality}
As it was already mentioned, 
the Kauffman polynomial obtained with
the parameters $(\a,s)$ and $(\a,-s^{-1})$
are equal. The corresponding BMW categories are also equal.
 From this we get 
 an
isomorphism between the constructed pre-modular
 categories. The image
of a  simple object $\l$
is $\l^\vee$. In fact the categories are equal; only the labelling of simple objects has changed.
 This  provides the ``level-rank'' duality isomorphism 
\begin{verse}
between $C^{n,k}$ and $C^{k,n}$, $B^{n,-k}$ and $B^{-k,n}$, 
$D^{n,k}$ and $D^{k,n}$;\\
between $C\!B^{n,k}$ and  $BC^{k,n}$, $B\!D^{n,k}$ and $D\!B^{k,n}$,
$C\!B^{n,-k}$ and $BC^{-k,n}$,\\ $B\!D^{-n,k}$ and $D\!B^{k,-n}$.
\end{verse}
\v8
Here we use that 
$\Gamma(D\!B^{k,n})=\{\l;\l_1+\l_2\leq 2n+1;\l_1^\vee+\l_2^\vee\leq 2k\}\ $.
In conclusion, up to the level-rank duality,
we  have obtained 
the following seven  series of pre-modular categories.

\begin{thm}
For $n,k\geq 1$
the categories $C^{n,k}$, $C\!B^{n,k}$, $C\!B^{n,-k}$
 $B^{n,-k}$, $B\!D^{n,k}$, $BD^{-n,k}$ and $D^{n,k}$ 
are pre-modular.
\end{thm}

 \section{  Modularization of the completed BMW categories}
\label{four}

In this section we discuss the modularization question for our series of pre-modular categories.

\subsection{Transparent simple objects} Let us first
note that  $\la \omega\ra=\sum_{\m\in {\Gamma (A)}}{\la\m\ra}^2$
 is nonzero
if $A$ is one of the pre-modular categories constructed in Section
\ref{trois};
the values of $\la \omega\ra$ are calculated e.g. in \cite{E}.
Therefore, the results of Section \ref{onetwo}
can be applied.

\begin{lem}\label{trans}
i) There is no non-trivial 
  transparent  simple object in the category $C^{n,k}$.

ii)  The non-trivial  transparent   simple objects are 
 $1^{2n+1},\, 2k+1,\, 1^{2n+1}\otimes (2k+1)$ in 
 $B^{n,-k}$ category;
$2k$, 
$1^{2n}$, $1^{2n}\otimes 2k$ in
$D^{n,k}$ category; 
 $2k$, $1^{2n+1}$, $1^{2n+1}\otimes
2k$ in
$B\!D^{n,k}$ and $B\!D^{-n,k}$ categories;
$2k+1$ in $C\!B^{n,k}$ and $C\!B^{n,-k}$ categories.


The quantum dimensions of these objects are equal to one.

\end{lem}

\begin{cor}\label{modcr} The   category $C^{n,k}$  with 
 $\Gamma(C^{n,k})$ as a representative set of simple objects
is modular.
\end{cor}

\noindent {\it Proof of the Lemma.} Recall that
a simple object $\l$ is transparent if and only
if for any (non-negligible) $\mu$ in the branching formula for $\l$,
the braiding coefficient 
is equal to one. 
Then $i)$ follows.

For $ii)$
we   verify that for each $\l$ mentioned in the lemma
all braiding coefficients  are  equal to one.
 Let us do it  in the $B^{n,-k}$ category for $\l=2k+1$.
Then  only $\m=2k$ appears in the branching formula for $\l$. We have
   $\l-\m=c$, $cn(c)=2k$ 
and  the braiding coefficient
is  $\a^{-2}s^{-4k}=s^{-4n-4k}=1$.
 Other cases can be done analogously. We see  that there
is no other transparent simple object in these categories.

The quantum dimensions can be calculated directly using (\ref{col}) and
$$\la j \ra=\frac{[0]_\a [1]_\a...[j-2]_\a([j-1]_\a +[j])}{[j]!}\; .$$
$\hfill\Box$

\begin{lem} For pre-modular categories constructed
 in Section \ref{trois}  the transparent simple objects 
form a group
under tensor multiplication. 
This group is isomorphic to $\mathbb Z_2\times \mathbb Z_2$
for  $D$, $B$, and $B\!D$ series
and to $\mathbb Z_2$ for  $C\!B$ series.

\end{lem}
\begin{proof}
It is sufficient to show that the  transparent
 simple objects have order 2, i.e.   any non-trivial  transparent  simple object
$t$ satisfies the equation: $t\otimes t\approx {\rm  trivial\,\, object}$.
Clearly, $t\otimes t$ contains the trivial object and
decomposes into a sum of transparent simple ones.
Comparing the quantum dimensions on the left and right hand side
of this decomposition formula we get the result.
\end{proof}

\v8

The twist coefficients  of  
the transparent objects listed in Lemma \ref{trans}
  are equal to $1$, except for the objects
$(2k+1)$ and $1^{2n+1}\otimes (2k+1)$
 in the  $B^{n,-k}$ category, 
$1^{2n+1}$ and $1^{2n+1}\otimes 2k$ in $B\!D^{-n,k}$ category,
and $(2k+1)$ in $C\!B^{n,-k}$ category,
whose twist coefficients are  $(-1)$.
Applying Brugui\`eres' criterion, we conclude.
\begin{cor}
The categories $D^{n,k}$, $B\!D^{n,k}$,
$C\!B^{n,k}$ are modularizable and $B^{n,-k}$, $B\!D^{-n,k}$, $C\!B^{n,-k}$
 are not modularizable.
\end{cor}

\noindent{\bf Remark.} The non-modularizable categories 
provide invariants of closed
 framed 3-manifolds (see \cite{sav}). 
Here a framing is a trivialization of the tangent bundle up to
isotopy. A choice of a framing is equivalent to the choice of
a spin structure and a 2-framing (or $p_1$-structure) on the 3-manifold.

\v8

\subsection{Modular categories 
 ${\bf \widetilde{ CB}}^{n,k}$,  ${\bf \widetilde{ BD}}^{n,k}$
 and  ${\bf \widetilde{ D}}^{n,k}$.}
Applying the modularization procedure described in Section \ref{one} to
the category $ C\!B^{n,k}$ we get 
the modular category $ \widetilde{C\!B}^{n,k}$ with the following
representative set of 
simple objects
$$ 
 \Gamma(\widetilde{ C\!B}^{n,k})=\{\l; \l_1 \leq k,
\l^\vee_1\leq n\}.$$ 
The stabilizer subgroup  for all
elements of $\Gamma({C\!B}^{n,k})$ is here trivial.

In the ${ B\!D}^{n,k}$ case, a simple object $\l$ with $\l_1=k$
has $Stab(\l)={\mathbb Z}_2$.
The algebra $End_{\widetilde{B\!D}^{n,k}}(\l)$ is two-dimensional.
It is generated by the tangle $a_\l$ having one free vertex 
colored by $2k$. We normalize it  such that $a^2_\l=1_\l.$
The minimal idempotents of $End_{\widetilde{B\!D}^{n,k}}(\l)$ are
$p^{\pm}_\l=1/2(\one\!_\l \pm a_\l)$. We define the simple objects
$\l_\pm$ 
  by means of  idempotents
$\tilde y_\m p^\pm_\l$. Their quantum dimensions are $\la\l\ra/2$.
As a result, 
$$
\Gamma(\widetilde{B\!D}^{n,k})
=\{\l;\l_1<k,\l^\vee_1\leq n\}\cup\{\l_\pm;\l_1=k, \l^\vee_1\leq n\}$$
is the representative set of simple objects
for the modular category $\widetilde{B\!D}^{n,k}$.

In the $D^{n,k}$ case, the diagrams belonging to the set $\Gamma_1=\{\l;
\l_1<k, \l^\vee_1<n\}$ have the  trivial stabilizer.
An object $\l$ from $\Gamma_2=\{\l;\l_1=k,\l^\vee_1<n \;\;\bigwedge\;\;
 \l_1<k,\l^\vee_1=n\}$
has the stabilizer equal to ${\mathbb Z}_2$.
We decompose it into $\l_\pm$ analogously to   the previous case.
An object  from $\Gamma_3=\{\l;\l_1=k,\l^\vee_1=n\}$
has the stabilizer of order 4.
The algebra $End_{\widetilde{D}^{n,k}}(\l)$, $\l\in \Gamma_3$,
 is 
either abelian or isomorphic
to the algebra of $2\times 2$
matrices. 

In the first case, $\l$ will decompose into the direct sum of
 four non-isomorphic
simple objects in the modular category $\widetilde {D}^{n,k}$.
 In the second case $\l$ will decompose into two isomorphic simple 
objects in $\widetilde {D}^{n,k}$.
It is a nontrivial open problem to decide which alternative holds for a given
$\l$. The answer may differ for distinct $\l$.
  To any  $\l\in \Gamma_3$ correspond $m_\l\in\{1,4\}$ simple
objects in $\Gamma({\widetilde D}^{n,k})$.
 If $m_\l=1$, we denote the object by $\hat\l$;
if $m_\l=4$, we denote the objects   
by ${}_\pm \l_\pm$.
Finally,  the representative set of simple objects
$\Gamma({\widetilde D}^{n,k})$  of the modular category
${\widetilde D}^{n,k}$ is
$$D_1=\Gamma_1 \cup \{\l_\pm;\l\in\Gamma_2\}\cup\{{}_\pm
\l_\pm; \l \in \Gamma_3,\ m_\l=4\}\cup\{
\hat\l; \l \in \Gamma_3,\ m_\l=1\}\ .$$

\section{ Verlinde formulas}
Recall that by  Turaev's work   any modular category  $\widetilde A$
with a set  $\Gamma$ of simple objects  gives rise to  a TQFT. 
The dimension of the TQFT module  associated with
 a genus $g$ closed surface is given  
by the Verlinde formula:
\be\label{ver}
d_g({\widetilde A})= \left(\sum_{\l\in \Gamma} {\la \l\ra}^2 \right)^{g-1}
\sum_{\l\in \Gamma} {\la \l\ra}^{2(1-g)} .\ee
In this section we calculate the dimensions of  TQFT modules
arising from the modular categories constructed above.

Let us introduce the notation ${[n]}_s=s^n-s^{-n}$ for $n\in \mathbb Z$.
\begin{thm}
i) The genus $g$ Verlinde formulas are 
$$d_g (C^{n,k})= 
(-(2n+2k+2))^{n(g-1)}\times$$
$$\times \sum_{n+k \geq l_1>...>l_n>0}
 \left(
{\prod^n_{j=1}{ {[2l_j]}_s} \prod_{1\leq i<j\leq n}
{{[l_i+l_j]}_s} {{[l_i-l_j    ]}_s}}\right)^{2(1-g)}; $$
\v8\v8

$$d_g (\widetilde{C\!B}^{n,k})=(-(2n+2k+1))^{n(g-1)}\times$$
$$\times\sum_{n+k\geq l_1>...>l_n>0}
 \left(
{\prod^n_{j=1}{ {[2l_j]}_s} \prod_{1\leq i<j\leq n}
{{[l_i+l_j]}_s} {{[l_i-l_j    ]}_s}}\right)^{2(1-g)}; $$



\v8\v8
$$\frac{d_g (\widetilde{B\!D}^{n,k})}{
(2n+2k-1)^{k(g-1)}}=
2\sum_{n+k-1\geq \a_1>...> \a_k>0}
 \left(
 \prod_{1\leq i<j\leq k}
{{[\a_i+\a_j]}_{ s}} {{[\a_i-\a_j    ]}_{ s}}\right)^{2(1-g)
}+$$ 
$$+ \sum_{n+k-1\geq \a_1>...> \a_k=0}
 \left(
 \prod_{1\leq i<j\leq k}
{{[\a_i+\a_j]}_{ s}} {{[\a_i-\a_j    ]}_{ s}}\right)^{2(1-g)}
.$$

 $$\frac{d_g ({\widetilde D}^{n,k})}{
{{(2n+2k-2)}}^{n(g-1)}}= 
\sum_{n+k-2\geq l_1>...> l_n= 0}
 \prod_{1\leq i<j\leq n}
\left({{[l_i+l_j]}_s} {{[l_i-l_j    ]}_s}\right)^{2(1-g)}
+ $$
$$+\;\;
2^{2g-1}\sum_{n+k-1=l_1>...> l_n= 0}
 \prod_{1\leq i<j\leq n}\left(
{{[l_i+l_j]}_s} {{[l_i-l_j    ]}_s}\right)^{2(1-g)}
+$$
$$ +\;\; 2 \sum_{n+k-2\geq l_1>...> l_n> 0}
 \prod_{1\leq i<j\leq n}\left(
{{[l_i+l_j]}_s} {{[l_i-l_j    ]}_s}\right)^{2(1-g)}
+$$
 $$+\;\;  \sum_{n+k-1=l_1>...> l_n> 0}(m_{(l-\delta)})^g
 \prod_{1\leq i<j\leq n}\left(
{{[l_i+l_j]}_s} {{[l_i-l_j    ]}_s}\right)^{2(1-g)}
. $$
Here $\delta=(n,n-1,\dots,1)$.
\v8\v8

ii) We have the following level-rank duality formulas.
$$d_g (C^{n,k})= d_g(C^{k,n})\;\;\;\;\;\;\;d_g ({\widetilde D}^{n,k})
= d_g({\widetilde D}^{k,n})$$
\end{thm}

\noindent {\bf Remark}.  The Verlinde formula for $C^{1,2}$
 calculates the number of the spin 
structures with Arf invariant zero  on the surface of genus $g$:
$d_g(C^{1,2})=2^{g-1}(1+2^{g})$. This fact should be interpreted 
via the corresponding TQFT, which is the is the one
 associated with the well known Rochlin invariant of spin $3$-manifolds
\cite{BM}.



\begin{proof} 
{\it i)} We substitute 
Propositions \ref{qdimC},
\ref{qdimB}, \ref{qdimD}  and calculations
of Sections 4.4-4.6 in \cite{E} into  (\ref{ver}).

Let us consider the $C^{n,k}$ case in  details. Here $\a=s^{-2k-1}$.
By Proposition \ref{qdimB}  and
 the calculations of  Section 4.5 in \cite{E}
we have
$$\sum_{\l\in \Gamma(C^{n,k})}{\la \l\ra}^2=
\frac{(-(2n+2k+2))^n}{\left(\prod^n_{j=1}[2n+2-2j]_s\prod_{1\leq i<j\leq n}[2n+2-i-j]_s
[j-i]_s\right)^2}\; .$$
Furthermore,
$$\sum_{\l\in \Gamma(C^{n,k})}{\la \l\ra}^{2(1-g)}=
\frac{\underset{n+k \geq l_1>...>l_n>0}{\sum}\left(
{\prod^n_{j=1}{ {[2l_j]}_s} \prod_{1\leq i<j\leq n}
{{[l_i+l_j]}_s} {{[l_i-l_j    ]}_s}}\right)^{2(1-g)}}
{\left(\prod^n_{j=1}[2n+2-2j]_s\underset{1\leq i<j\leq n}{\prod}[2n+2-i-j]_s
[j-i]_s\right)^{2(1-g)}}\; .$$
Here we used  the bijection 
$\Gamma(C^{n,k})\to T:=\{(l_1,...,l_n), n+k \geq l_1>...>l_n>0\}$
 sending $\l$ to $\l+(n,n-1,...,1)$.
Substituting  the last two formulas
 into (\ref{ver}) we get the result.

 For the third formula we use that

$$\sum_{\l\in \Gamma({B\!D}^{n,k})}{\la \l\ra}^2=
4 \sum_{\l\in \Gamma(\widetilde{BD}^{n,k})}{\la \l\ra}^2 .$$
This is because  the action of the group 
${\mathbb Z}_2 \times {\mathbb Z}_2$ of the transparent objects on 
$\{\l; \l_1<k,\l^\vee_1\leq n\}$ 
preserves the quantum dimension 
 and
$\la\l_\pm\ra=1/2 \la \l\ra$.

{\it ii)} 
 By (\ref{sym}) we have for any $p\in {\mathbb{ Z}}$
$$\underset{\l^\vee_1\leq n}{\sum_{\l_1\leq k}}{\la \l\ra}^p_{-s^{2n+1},s}=
\underset{\l^\vee_1\leq n}{\sum_{\l_1\leq k}}{\la \l^\vee\ra}^p_{-s^{2n+1},
-s^{-1}}=
\underset{\l_1\leq n}{\sum_{\l^\vee_1\leq k}}{\la \l\ra}^p_{- s^{2k+1},
 s}\; .$$
The second formula can be shown analogously.

 \end{proof}

\section{Refinements} 

In this section we construct spin and cohomological refinements
of the quantum invariants arising from the  modular category
$C^{n,k}$.

We work here
in $C^{n}$-$C^{k}$ 
 specialization,
 i.e.  $\a=s^{-2k-1}$, $s$ is a primitive $2l$th root of
unity, $l=2n+2k+2$. Recall
$\Gamma(C^{n,k})=\{\l; \l_1\leq k, \l^\vee_1\leq n\}\, .$
Let us introduce
 a ${\mathbb{Z}}_2$-grading on 
the category  $C^{n,k}$ 
corresponding to the parity of the number of cells in  
Young diagrams. According to this grading,
 we decompose  the Kirby element: $\omega=\omega_0 + \omega_1$.
\begin{lem}\label{71}
Let  $U_\varepsilon(\l)$ be the $\varepsilon$-framed unknot
colored with $\l$ and $\varepsilon =\pm 1$.

i) For  $kn=2\mod 4$, we have
$\la U_\varepsilon(\omega_0)\ra=0$.

ii)  For  $kn=0\mod 4$, we have 
$\la U_\varepsilon(\omega_1)\ra=0$.
\end{lem}

\begin{proof}
Let us  call the {\it graded sliding property} the equality  
drawn in Proposition \ref{slidprop}
 by replacing $\omega$ on the left-hand side by $\omega_\n$
and $\omega$ on the right-hand side by $\omega_{\n+1}$ with  $\n=0,1$.
The proof of this identity can be adapted from the one of this proposition.

Using twice  the graded
 sliding property, we can see that the morphism drawn below
is nonzero only if $\l=0$ or $\l=k^n$.
$$\psdiag{8}{24}{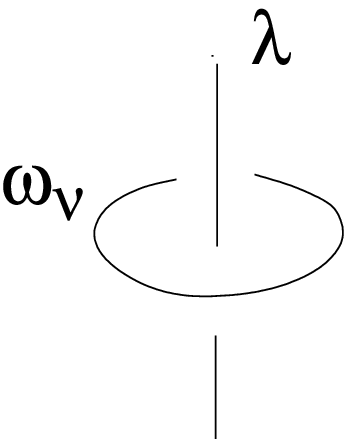}$$
Then
\be\label{43}
\la U_1(\omega_\n)\ra \la U_{-1}(\omega_\n)\ra=
\la H_{1,0}(\omega_0,\omega_\n)\ra=(1+\alpha^{kn}s^{nk(k-n)}s^{l\n})\la 
\omega_\n\ra\; \ee
where  $H_{1,0}(\omega_0,\omega_\n)$
is the Hopf link whose $\omega_0$-colored component has framing 1
and $\omega_\n$-colored one is $0$-framed.
The first equality is due to
the graded sliding property. In the second one we use the twist
and braiding coefficients for $\l=0,k^n$ and the fact that
$$\sum_{c\in k^n}cn(c)=\frac{nk}{2}(k-n)\; .$$ 
 Substituting the values of $\alpha$ and $l$ 
into (\ref{43})
we get the result.
\end{proof}

The following statement is the direct consequence of  this lemma and 
the construction of refined invariants described in \cite[Section 4]{Bl}.
\begin{thm}\label{72} The
quantum invariants arising from the modular category $C^{n,k}$ 
 can be written as  sums  of  refined invariants corresponding to 
different 
spin structures if $kn=2\mod 4$ and to  $\mathbb{Z}_2$-cohomology classes
if $kn=0\mod 4$.
\end{thm}

One can show by the same method 
that other categories do not provide refined invariants. 

\section{Comparison with the quantum group approach}

The aim of this section is to find a correspondence between
  pre-modular categories that have been 
constructed in Section \ref{trois} and those that arise from
the quantum group method.

\subsection{Modular categories from quantum groups}

We keep  notation of \cite{Le1}, \cite{ Le2}.
Let $(a_{ij})_{1\leq i,j\leq l} $ be the Cartan matrix of a simple
complex Lie algebra $\mathfrak g$.  
There are relatively prime integers $d_1,...,d_l$ in $\{1,2,3\}$ such that
the matrix $(d_i a_{ij})$ is symmetric. Let $d=max(d_i)$.
We fix a Cartan subalgebra $\mathfrak h$ of $\mathfrak g$ and
fundamental roots  $\a_1, \a_2,...,\a_l$
in the dual space ${\mathfrak h}^*$. 
Let ${\mathfrak h}^*_{\mathbb R}$ be the $\mathbb R$-vector space 
spanned by the fundamental roots. The root lattice $Y$ is the 
$\mathbb Z$-lattice generated by $\a_i$, $i=1,...,l$.
We define an inner product on ${\mathfrak h}^\ast_{\mathbb R}$ 
by  $(\a_i|\a_j)=
d_i a_{ij}$. Then $(\a|\a)=2$ for every short root $\a$.
The inner product normalized such that every long root has length two
will be denoted by $(.|.)'$.  We have $(.|.)'=(.|.)/d$. 
Let $\l_1, ...,\l_l$ be the fundamental weights, 
then $(\l_i|\a_j)=d_i\delta_{ij}$.
The weight lattice $X$ is the $\mathbb Z$-lattice generated by $\l_1,...,\l_l$.
Let $\rho=\l_1+...+\l_l$.
The Weyl chamber is
defined by
$C=\{x\in {\mathfrak h}^*_{\mathbb R}; (x|\a_i)\geq 0, i=1,...,l\}\, .$
Let us denote by $\a_0$ (resp. $\beta_0$) the short 
(resp. the long) root in the Weyl chamber $C$.

Let $U_q({\mathfrak g})$ be the quantum group associated with $\mathfrak g$
and  $q$ be  a primitive root of unity of order 
$r$ (notation coincides with \cite[Section 1]{Le2}).
Let  $h^\vee$ be the dual Coxeter number.
The case when $r\geq d h^\vee$ is divisible by $d$ was
mainly studied in the literature. In that case, 
simple  $U_q({\mathfrak g})$-modules 
corresponding to weights in
$$C_L=\{x\in C; (x|\beta_0)'\leq L\}$$
form a pre-modular category \cite{Kir}.
Here $L:=r/d-h^\vee$ is the level of the category.
The quantum dimension of  $\m\in  X$ is given by
\be\label{qqqddd}
\la\m\ra=\prod_{\text{\tiny positive roots}\; \a} \frac{v^{(\m+\rho|\a)}-
v^{-(\m+\rho|\a)}}{v^{(\rho|\a)}-v^{-(\rho|\a)}},\ee
its twist coefficient is $v^{(\m+2\rho|\m)}$, where $v^2=q$.
The modularization of these categories was studied in \cite{sav1}.

In the case when $(r,d)=1$ and $r>h$ ($h$ is the Coxeter number), 
pre-modular categories
can also be constructed \cite{Le1}.
The set of simple objects corresponds to weights
in
$$C'_L=\{x\in C; (x|\a_0)\leq L\}$$
with $L:=r-h$. 
Le showed that if $(r, d\;{\rm  det}(a_{ij}))=1$ the set of modules
in $C'_L\cap Y$ generates a modular category.

We say that
two pre-modular categories are equivalent  if there exists
a bijection between their sets of simple objects providing an equality
of the corresponding colored link invariants.
For modularizable categories 
this implies that  the associated  TQFT's are isomorphic
 (compare \cite[III,3.3]{Tu}).

\subsection{Comparison of C cases} 
Any weight $\m\in C$ of $C_n$ is of the form
$\m=\l_1 e_1 +...+\l_n e_n$ with $(e_i|e_j)=\delta_{ij}$
and integers $\l_1\geq...\geq\l_n\geq 0$ (compare \cite[p.293]{On}).
 With any $\m$ a Young diagram
$\l=(\l_1,...,\l_n)$ can be associated.

\begin{thm} The pre-modular categories associated with $U_q(C_n)$ 
and $U_q(C_k)$ at a primitive root of unity 
$q$ of order $r=2n+2k+2$ are equivalent to $C^{n,k}$.
\end{thm}
\begin{proof}

A colored $m$-component link invariant of a pre-modular category
$A$ with $\Gamma(A)$ as a representative set of simple objects
can be considered as a multilinear function from ${\Gamma(A)}^m$ to $\fk $.
Here we supply $\Gamma(A)$ with a ring structure by considering
direct sums and tensor products. It is easy to see from the previous 
discussion, that there exists
an isomorphism  between such  rings in our case. Indeed,
for $U_q(C_n)$ we have $d=2$, $h^\vee=n+1$, $L=r/2-n-1=k$,
 $\beta_0=2e_1$,  
and $C_L\cap X=\{\m; \l_1\leq k\}$. After the 
identification of $\m$ with $\l$,
this coincides with
$\Gamma(C^{n,k})$. The ring structure is preserved under this identification.
Furthermore, it is known that these rings are generated by
 the fundamental module corresponding to $\m=e_1$ and the object  $\l=1$.
Therefore, it is sufficient to verify the equality of invariants
colored by these two objects.
The fact, that the link invariant associated 
with this fundamental module
is a specialization of the Kauffman polynomial  was shown in
\cite{Tu1}. In order to identify the parameters, compare the 
quantum dimensions of simple objects given by 
(\ref{qqqddd}) and Proposition \ref{qdimC}. We show that $s^2=q$.
The equivalence between  $U_q(C_k)$ 
and $C^{k,n}$ can be shown analogously. Then we use
the level-rank duality.
\end{proof}

Analogously, the category
$C\!B^{n,k}$  is equivalent to $U_q(C_n)$ with
$r=l=2n+2k+1$. Indeed, we have  $(r,d)=1$, $h=2n$,
$L=r-2n=2k+1$, $\a_0=e_1+e_2$, and
$C'_L\cap X=\{\m;\l_1+\l_2\leq 2k+1\}$.

\subsection{Comparison of  B cases}

 Any weight $\m\in C$ of $B_n$ can be written in the form
$\m=\l_1 e_1 +...+\l_n e_n$, where $(e_i|e_j)=2 \delta_{ij}$
 and half-integers $\l_1\geq...\geq \l_n\geq 0$.
 If   $\l_i\in\mathbb N$,
$i=1,...,n$, the partition
 $\l=(\l_1,...,\l_n)$ 
defines a Young diagram associated  with $\m$.
If at least one of $\l_i$ is a half-integer, we call $\m$ a spin module.
Our construction of simple objects
can be considered as  a quantum
analog of the Weyl construction and
it does not produce spin modules.
Let us compare the quantum dimension and/or twist coefficient
of a non-spin module $\m$ and the corresponding Young diagram $\l$
given by  (\ref{qqqddd}) and Proposition \ref{qdimB}.
They coincide if  $v^2=q=s$. 

Let us first consider the case when $r$ is even and $r>4n$. 
Here $h^\vee=2n-1$.
Let $r=4n+4k$ with $k\geq 1$. We have $\beta_0=e_1+e_2$
and  $L=r/2-2n+1$. Then
$$C_L\cap X=\{\m; \l_1+\l_2\leq 2k+1 \}.$$
We conclude that the quotient by spin modules of the pre-modular category
for $B_n$ at $(4n+4k)$th root of unity is equivalent to
 the modular extension of $B^{n,-k}$ by $G'$ generated by $1^{2n+1}$.
Using the level-rank duality, we get the equivalence of 
the pre-modular category for $B_k$ at $(4n+4k)$th 
root of unity with  the modular extension of $B^{n,-k}$ by $G''$ 
generated by ${2k+1}$.

Let us put $r=4n+4k-2$ with $k\geq 1$. Then we get analogously that
 the category $B\!D^{n,k}$ is equivalent
to the quotient (by spin modules)
of the pre-modular category for $B_n$ at $(4n+4k-2)$th  root 
of unity.

For odd $r>h=2n$ we set $r=2n+2k+1$ with $k\geq  1$.
Then  $(r,d)=1$, $\a_0=e_1$ and
$L=r-2n=2k+1$. We have $C'_L\cap X=\{\m;\l_1\leq k+1/2\}$.
 We see that the quotient by spin modules of the pre-modular category
for $B_n$ on $(2n+2k+1)$th root of unity 
is equivalent to  $\widetilde{ { BC}}^{n,k}$.

\subsection{Comparison of D cases}

Any weight of $D_n$ can be written in the form
$\m_\pm=\l_1 e_1 +...+\l_{n-1} e_{n-1} \pm \l_n e_n $ 
with $(e_i|e_j)=\delta_{ij}$
and half-integers $\l_1\geq...\geq\l_n\geq 0$.
Here we have $v=s$, $d=1$, $h=2n-2$ and $\beta_0=e_1+e_2$.
Setting $r=2n+2k-2$ with $k\geq 1$
we get
$$C_L\cap X= \{\m_\pm;\l_1+\l_2\leq 2k \}$$
For non spin modules, this coincides with 
the set of simple objects of the modular extension of $D^{n,k}$ by
$G'$ generated by $1^{2n}$.
Therefore, the modular categories for $D_n$ and $D_k$
at $(2n+2k-2)$th root of unity are equivalent to
$\widetilde D^{n,k}$.

For $r=2n+2k-1$ ($k\geq 1)$ we get
 that the modular extension of ${ D\!B}^{n,k}$ (isomorphic to
 $B\!D^{k,n}$)  by $G'$ as above
is equivalent to the quotient by spin modules
 of the  pre-modular category for $D_n$. 


\v8
As a result, any  pre-modular category defined in Section \ref{trois}
 is equivalent to a  
quantum group category.
Moreover, our  categories
produce a complete set of invariants
that  can be obtained
from  quantum groups of types B,C and D  by using non-spin modules.

\bibliographystyle{amsplain}

\end{document}